\documentclass[preprint]{elsarticle}

\usepackage[english]{babel}
\usepackage{amssymb}
\usepackage{subfigure}
\usepackage[active]{srcltx}     
\usepackage[dvips]{epsfig}
\usepackage[english]{babel}
\usepackage{graphics}
\usepackage{graphicx}
\usepackage{epsfig}
\usepackage{amsmath,amssymb}
\usepackage{amsfonts}
\graphicspath{{./Figures/}}
\numberwithin{equation}{section}
\usepackage{color}
\usepackage{verbatim}
\usepackage{dcolumn}
\usepackage{changes}
\usepackage{placeins}
\usepackage{url}
\usepackage{multirow}
\usepackage[group-separator={,}]{siunitx}
\usepackage{bm}
\usepackage[mathscr]{euscript}
\usepackage{placeins}
\journal{Int. J. for Num. Methods in Engrg.}

\begin{document}

\begin{frontmatter}

\title{Virtual Element Method: an equilibrium-based \\stress recovery procedure}

\author{E. Artioli\corref{cor1}$^a$}
\author{S. de Miranda\corref{cor2}$^b$}
\author{C. Lovadina\corref{cor4}$^{c,d}$}
\author{L. Patruno $^b$ \corref{cor33}}
\ead{luca.patruno@unibo.it}

\cortext[cor33]{Corresponding author}

\address{$^a$Department of Civil Engineering and Computer Science, University of Rome Tor Vergata, Via del Politecnico 1, 00133 Rome, Italy}
\address{$^b$DICAM, University of Bologna, Viale Risorgimento 2, 40136 Bologna, Italy}
\address{$^c$Dipartimento di Matematica, Universit\`a degli Studi di Milano, Via Saldini 50, 20133 Milano, Italy}
\address{$^d$IMATI del CNR, Via Ferrata 1, 27100 Pavia, Italy}

\begin{abstract}
Within the framework of the displacement-based Virtual Element Method (VEM) for plane elasticity a significant problem is represented by an accurate evaluation of the stress field. In particular, in the classical VEM formulation, a suitable operator which maps to the strain field is introduced in order to allow the calculation of the stiffness matrix. The stress field is then computed using that strian field, by using the constitutive law. Considering for example a first-order formulation for a homogeneous material, strains are locally mapped onto constant functions, and stresses are accordingly piecewise constant. However, the virtual displacements might engender more complex strain fields for polygons which are not triangles. In this paper, Recovery by Compatibility in Patches is used in order to mitigate such an effect and, thus, enhance the accuracy of the recovered stress field. The procedure is simple, efficient and can be readily implemented in existing codes. Numerical tests confirm the soundness of the proposed approach.
\end{abstract}


\end{frontmatter}

\section{Introduction} \label{sec:introduction}
The Virtual Element Method (VEM) is a relatively new and very powerful discretisation scheme which is rapidly attracting the interest of the scientific community. The technique is well-known for its formal elegance and flexibility, which allows to adopt meshes composed of general polygons/polyhedra as well as allowing the presence of hanging nodes and nonconforming grids \cite{volley,Hitchhikers,BeiraodaVeiga-Brezzi-Marini:2013,nonconforming}.
The method, originally proposed  in 2012 in its displacement-based formulation and presented for the Laplace operator in a two-dimensional context \cite{volley}, is rapidly developing and has been already applied to numerous physical problems \cite{BeiraodaVeiga-Brezzi-Marini:2013,Antonietti-BeiraodaVeiga-Mora-Verani:20XX,Brezzi-Marini:2012,Artioli2017a,Artioli2017b,WRRH,BCP,TPPM10,ABSV2017,Berrone-VEM,Helmholtz-VEM,Steklov-VEM} as well as extended to three-dimensional cases \cite{Paulino-VEM} and mixed formulations \cite{Brezzi-Falk-Marini,Artioli2017HR}.

Differently from the finite element framework, in the virtual element technique the local polynomial approximation concept is abandoned and the functions used for the discretization procedure are not known pointwise, in general. However, a careful selection of the degrees of freedom and assuming that the unknown fields satisfy appropriate differential equations, make it possible is to compute the stiffness matrix and the load vector.

In particular, for plane elasticity problems, classical VEM formulations assume that displacements are completely known only at the element boundaries while the internal field is unknown and, for high-order schemes, characterised only by means of internal degrees of freedom representing appropriate integral quantities.

Due to the virtual nature of the displacement field, strains (the symmetric gradient of the displacements) are not directly computable and need to be further approximated, in order to compute the stiffness matrix. The standard procedure consists in defining a strain field, polynomial inside the element, which can be determined starting from the displacement degrees od freedom.
For example, in the case of the first-order formulation, the aforementioned strategy leads to constant strains inside the elements, irrespective of the number of element edges, and thus irrespective of the number of displacement degrees of freedom. Therefore, applying the constitutive law, the stress field is essentially piecewise constant as well. Especially for meshes consisting of polygons with numerous edges, this often leads to not completely satisfactory results.

In this paper, a modified version of the Recovery by Compatibility in Patches (RCP) \cite{Ubertini2004,Ubertini2006} is used to devise an alternative procedure, which allows to partly alleviate the above drawback and compute accurate stresses in the displacement-based VEM framework. Patch-based stress recovery techniques, initially inspired by the Superconvergent Patch Recovery (SPR) procedure proposed by Zienkiewicz and Zhu \cite{Zienkiewicz1992a,Zienkiewicz1992b}, have quickly developed and found numerous applications, ranging from mesh adaptivity \cite{Boroomand1999b,Castellazzi2010} to the recovery of accurate stresses for crack propagation in XFEM analyses \cite{Prange2012}.
In contrast with classical approaches, the common characteristic of such techniques is to recover enhanced stress fields by considering groups of surrounding elements (i.e. patches). In the original SPR procedure this is obtained by calculating stresses as a least-square interpolation of values evaluated at superconvergent points. Later, the procedure has been refined \cite{Blacker1994,Wiberg1997,Lee1997,Park1999,Kvamsdal1998,Okstad1999} introducing, for example, equilibrium constraints in order to lead to more accurate results. Alternative procedures, based on variational principles such as the Recovery by Equilibrium in Patches (REP) \cite{Boroomand1997a,Boroomand1997b,Boroomand1999a}, have been proved to obtain extremely interesting results.

RCP has been originally proposed as a patch-based recovery technique which, by enforcing equilibrium and relaxing compatibility, allows to enhance the stresses obtained from standard displacement-based finite elements schemes.

The key idea underlying RCP is to minimise the complementary energy associated to a patch of elements treated as an isolated system. On each patch, stresses are obtained as a linear combination of self-equilibrated stress modes enriched by an appropriate particular solution, and the explicit knowledge of the displacements is required only along the patch boundaries. Therefore, the RCP approach is naturally well-suited for the virtual element schemes.

In this paper, the application of RCP in the context of VEM is analysed by proposing two approaches. In the first one, a degenerate patch composed of a single polygon is considered. For triangles, this does not significantly improve the results obtained using the classical strain computation.  However, the benefit becomes apparent when the number of edges increases. In the second form, RCP is adopted by considering patches of elements in analogy to its original formulation developed in the context of finite elements schemes.

The paper is organised as follows. In Section \ref{sec:vem}, the VEM formulation for plane elasticity problems is briefly recalled. Then, the RCP procedure is introduced in Section \ref{sec:rcp}. Numerical results are presented for a wide selection of test cases in Section \ref{sec:num} and, finally, conclusions are drawn in Section \ref{sec:con}.

\section{Virtual Element Method}
\label{sec:vem}
Consider a body occupying a region $\Omega$ of the two-dimensional space on which a reference system $(O,x,y)$ is introduced. Let us denote with the symbol $\partial \Omega$ the boundary of such a body and consider the case of homogeneous Dirichlet boundary conditions. Other types of boundary conditions can be treated in the same spirit of finite element schemes.  Indicating as $\bm{u}(x,y) = [u,v]^T$, the displacement field, the associated strains are defined as

\begin{equation}
\bm{\varepsilon}(\bm{u}) = \bm{D}\bm{u}   \qquad \text{with} \qquad \bm{D}=\left[
                                                                          \begin{array}{cc}
                                                                            \partial_x & 0 \\
                                                                            0 & \partial_y \\
                                                                            \partial_y & \partial_x \\
                                                                          \end{array}
                                                                        \right],
\end{equation}

\noindent where $\bm{D}$ represents the differential compatibility operator. The plane elasticity problem can be thus expressed as: find $\bm{u} \in \bm{V}$ such that
\begin{equation}
\int_\Omega \bm{\varepsilon}(\bm{v})^T \bm{C} \bm{\varepsilon}(\bm{u})= \int_\Omega \bm{v}^T \bm{b} \qquad \forall \bm{v} \in \bm{V}, \label{eq:prob}
\end{equation}
\noindent where $\bm{C}$ is the elastic matrix, $\bm{V}$ is the space of the kinematically admissible displacements, $\bm{b}$ is the applied external load.

In order to build a VEM scheme, it is assumed that $\Omega$ is subdivided into non-overlapping polygons with straight edges. In the following, $E$ indicates one and each of such polygons, $\partial E$ denotes its boundary while $e$ denotes one and each of the edges of $E$. Moreover, we here assume a homogeneous material, but inhomogeneous bodies can be treated as well (cf. for instance \cite{Artioli2017HR, AdMLP_Ho}) We introduce a local approximation space, $\bm{V}_{h|E}$, whose elements are typically not known pointwise, but are determined by means of a suitable set of degrees of freedom (thus, $\bm{V}_{h|E}$ is said to be {\em virtual}).

In the following, reference is made to a first-order VEM discretisation scheme. The space $\bm{V}_{h|E}$ is chosen so that
\begin{equation}
\bm{V}_{h|E} = \left\{ \bm{v}_h \in H^1(E) \cap C^0(E): \mathbf{D^\star}(\bm{C} \bm{\varepsilon}(\bm{v}_h))=0,\; \bm{v}_{h|e} \in \bm{P}_{1}(e) \quad \forall e \in \partial E  \right\},
\label{eq:space}
\end{equation}

\noindent where $\mathbf{D^\star}$ is the formal adjoint of $\mathbf{D}$ while $\bm{P}_{1}(e)$ indicates the space of linear functions on the edges of $E$.

\begin{figure}[h!]
	\centering
	\includegraphics[width=0.9\textwidth,trim = 0mm 0mm 0mm 0mm, clip]{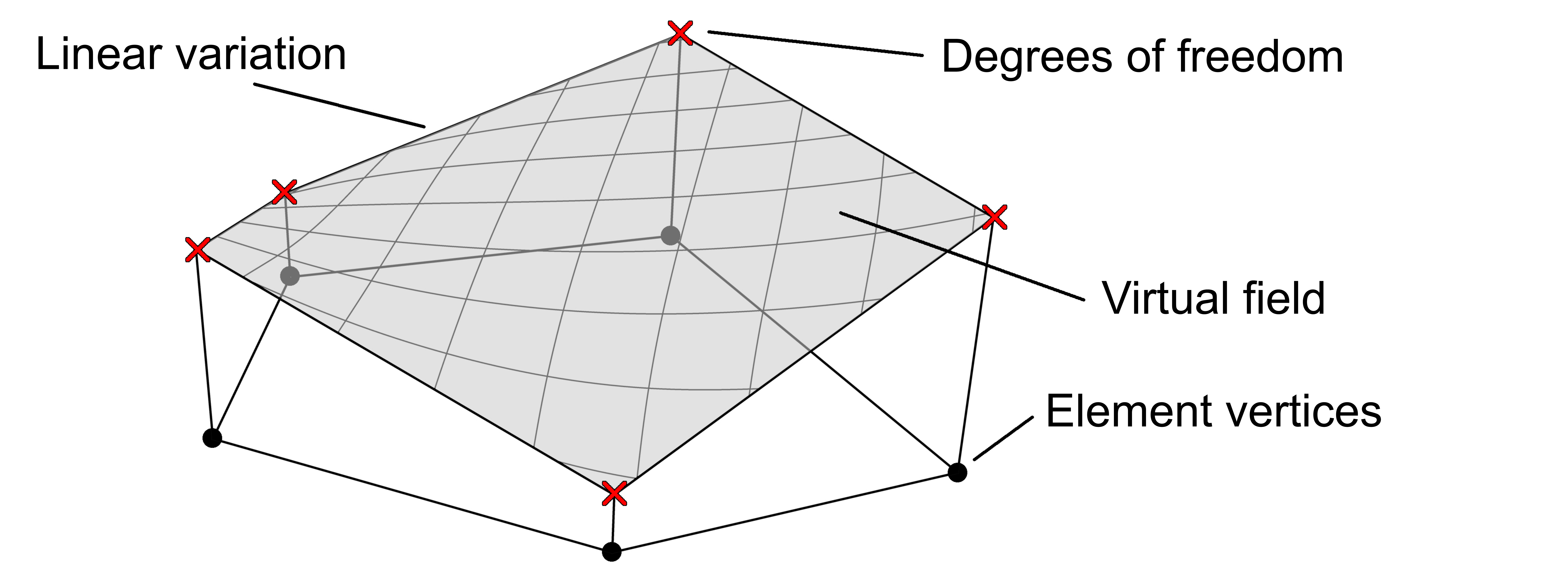}
	\caption{Sketch of the virtual function and of the considered degrees of freedom.}
	\label{fig:vField}
\end{figure}

A set of degrees of freedom which can be used in order to describe such a space is represented by the values of the displacements at the element vertices. It must be noticed that, thanks to the requirements expressed in Eq. \eqref{eq:space}, displacements must be linear on edges so that the vertex values completely characterise the displacement field on $\partial E$. A sketch of the virtual field and of the considered degrees of freedom is provided in Fig. \ref{fig:vField}.

Due to the virtual nature of the approximating functions, it is not possible to directly compute the terms of Eq. \eqref{eq:prob}, and especially the left-hand side. It is thus necessary to introduce an operator, $\bm{\Pi}: \bm{V}_{h|E}\to \bm{P}_0(E)$, defined as
\begin{equation}
\int_E \bm{\Pi}(\bm{v}_h)^T \bm{\varepsilon}^0 = \int_E \bm{\varepsilon}(\bm{v}_h)^T \bm{\varepsilon}^0 \qquad \forall \bm{\varepsilon}^0 \in \bm{P}_0(E).
\label{eq:proj}
\end{equation}
Above, $\bm{P}_0(E)$ denotes the space of constant fields over the element. It can be observed that $\bm{\Pi}$ is an operator which maps the virtual displacements onto constant strain fields. A typical displacement and strain field can be formally represented as
\begin{equation}
\bm{v}_h = \bm{N}^V \tilde{\bm{v}}_h, \qquad \bm{\varepsilon}^0 = \bm{N}^P \tilde{\bm{\varepsilon}},
\end{equation}
\noindent where the vector $\tilde{\bm{v}}_h$ collects the degrees of freedom (i.e. the values of the displacements at the polygon vertices) while $\tilde{\bm{\varepsilon}}$ is a vector collecting the degrees of freedom of the approximated strain field. The term $\bm{N}^V$ collects virtual displacement shape functions and, thus, is not explicitly known, while $\bm{N}^P$, still considering a first order scheme, takes the form
\begin{equation}
\bm{N}^P = \left[ \begin{array}{ccc}
                    1 & 0 & 0 \\
                    0 & 1 & 0 \\
                    0 & 0 & 1
                  \end{array} \right].
\end{equation}

Indicating as $\bm{\Pi}^m$ the $\bm{\Pi}$ operator expressed as a matrix, it is possible to write the discrete strains as
\begin{equation}
\bm{\Pi}(\bm{v}_h) = \bm{N}^P \bm{\Pi}^m \tilde{\bm{v}}_h. \label{eq:eproj}
\end{equation}

The operator $\bm{\Pi}^m$ can be thus evaluated by means of Eq. \eqref{eq:proj}. In particular, substituting Eq. \eqref{eq:eproj} into Eq. \eqref{eq:proj}, we obtain
\begin{equation}
\int_E (\bm{N}^P\bm{\Pi}^m\tilde{\bm{v}}_h)^T \bm{N}^P \tilde{\bm{\varepsilon}} = \int_E \left[\bm{\varepsilon}(\bm{N}^V \tilde{\bm{v}}_h)\right]^T \bm{N}^P \tilde{\bm{\varepsilon}} \qquad \forall \tilde{\bm{\varepsilon}} \in \mathbb{R}^3,
\end{equation}

\noindent which, integrated by parts, leads to
\begin{equation}
\int_E (\bm{N}^P\bm{\Pi}^m\tilde{\bm{v}}_h)^T \bm{N}^P \tilde{\bm{\varepsilon}} = \int_{\partial E}  (\bm{N}^V \tilde{\bm{v}}_h )^T \bm{N}_E\bm{N}^P\tilde{\bm{\varepsilon}} \qquad \forall \tilde{\bm{\varepsilon}} \in \mathbb{R}^3,
\label{eq:projeq}
\end{equation}

\noindent where, denoting the outward normal as $\bm{n}=[n_x,n_y]^T$, $\bm{N}_E$ is
\begin{equation}
\bm{N}_E = \left[ \begin{array}{ccc}
             n_x & 0 & n_y \\
             0 & n_y & n_x
           \end{array} \right].
\end{equation}

\noindent Equation \eqref{eq:projeq} can be rewritten as
\begin{equation}
\tilde{\bm{\varepsilon}}^T \mathcal{\bm{G}} \bm{\Pi}^m \tilde{\bm{v}}_h = \tilde{\bm{\varepsilon}}^T \mathcal{\bm{B}} \tilde{\bm{v}}_h  \qquad \forall \tilde{\bm{\varepsilon}} \in \mathbb{R}^3, \label{eq:pim}
\end{equation}
\noindent where
\begin{equation}
\mathcal{\bm{G}} = \int_E (\bm{N}^P)^T \bm{N}^P, \qquad \mathcal{\bm{B}}  = \int_{\partial E}  (\bm{N}_E \bm{N}^P)^T \bm{N}^V.
\end{equation}

Equation \eqref{eq:pim} allows to compute $\bm{\Pi}^m$ as the solution of a linear system while matrices $\mathcal{\bm{G}}$ and $\mathcal{\bm{B}}$ are computable because displacements appear only at the boundary where they are explicitly known. This leads to the evaluation of $\bm{\Pi}^m$ as
\begin{equation}
\bm{\Pi}^m = \mathcal{\bm{G}}^{-1} \mathcal{\bm{B}},
\end{equation}
\noindent and, thus, to approximate the l.h.s. of Eq. \eqref{eq:prob} as
\begin{equation}
\int_\Omega \bm{\varepsilon}(\bm{v})^T \bm{C} \bm{\varepsilon}(\bm{u}) \approx \int_\Omega (\bm{\Pi}^m{\tilde{\bm{v}}})^T \bm{C}  \bm{\Pi}^m{\tilde{\bm{u}}}.
\end{equation}

\noindent After algebraic manipulations, it is possible to identify the stiffness matrix, $\bm{K}_c$, as
\begin{equation}
\bm{K}_c = \mathcal{\bm{B}}^T \mathcal{\bm{G}}^{-T} \left[ \int_E (\bm{N}^P)^T \bm{C} \bm{N}^P \right] \mathcal{\bm{G}}^{-1} \mathcal{\bm{B}}.
\end{equation}

\noindent In order to ensure that $\bm{K}_c$ has the correct rank, a stabilisation term $\bm{K}_s$ (not here further discussed)  must be added to $\bm{K}_c$ so that the stiffness matrix is usually redefined as $\bm{K} = \bm{K}_c + \bm{K}_s$. For further details, we refer to \cite{Hitchhikers,Artioli2017a}, for example.

It is important to notice that the operator $\bm{\Pi}^m$ leads to constant strain fields while the virtual displacements are assumed to be linear over the element edges, and lead to divergence-free stress distributions (cf. Eq. \eqref{eq:space} and Fig. \ref{fig:vField}).

For triangular elements, the virtual (low-order) formulation is equivalent to standard (low-order) triangular finite elements. However, when polygons with a higher number of edges are considered, the degrees of freedom for the displacements could, in principle, give rise to stress/strain fields richer than the constants. A kind of loss of information is thus expected when using operator $\bm{\Pi}^m$ to compute the stresses.

In the following, RCP is introduced as a postprocessing technique for the VEM procedure, with the aim of improving the accuracy of the obtained stresses.

\section{Recovery by Compatibility in Patches}
\label{sec:rcp}
The key idea of RCP is to minimise the complementary energy associated to a patch of elements treated as an isolated system \cite{Ubertini2004}. Such operation allows to enforce equilibrium, which might be desirable in some applications (see for instance \cite{deMiranda2012}), as well as to increase the accuracy of the obtained stress field.
In this work, RCP is presented for degenerate patches composed of a single virtual element. The generalisation of the procedure to the case of patches composed of more than one element does not introduce substantial modifications.

The complementary energy associated to $E$ can be expressed as
\begin{equation}
J^c(\bm{\sigma}^\ast) = \frac{1}{2} \int_E \bm{\sigma}^{\ast T} \bm{C}^{-1} \bm{\sigma}^\ast - \int_{\partial E} \bm{u}^T \bm{N}_E \bm{\sigma}^\ast,
\label{eq:cener}
\end{equation}

\noindent where $\bm{\sigma}^\ast$ is the unknown stress field while $\bm{u}^T$ is the known displacement field (i.e. the displacement field computed by the displacement-based VEM analysis). It must be noticed that the displacement field appears only in the boundary term. In order to minimise $J^c$, the stress field is written as the sum of two contributions:
\begin{equation}
\bm{\sigma}^\ast = \bm{\sigma}^\ast_h + \bm{\sigma}^\ast_p,
\end{equation}
\noindent which correspond to its divergence-free part and an appropriate particular solution which guarantees equilibrium with the applied external loads. The homogeneous part $\bm{\sigma}^\ast_h$ is approximated as
\begin{equation}
\bm{\sigma}_h^\ast  = \bm{P} \bm{\beta},
\label{eq:stressModes}
\end{equation}
\noindent where $\bm{P}$ is a matrix of preselected self-equilibrated stress modes, while $\bm{\beta}$ is the vector of the unknown parameters. In particular, if a linear approximation is adopted for $\bm{\sigma}_h$, $\bm{P}$ takes the form
\begin{equation}
\bm{P} = \left[ \begin{array}{ccccccc}
                  1 & 0 & 0 & y & 0 & x & 0 \\
                  0 & 1 & 0 & 0 & x & 0 & y \\
                  0 & 0 & 1 & 0 & 0 & -y & -x
                \end{array} \right].
\label{eq:pmatr}
\end{equation}

The particular solutions, $\bm{\sigma}_p^\ast$, appearing in Eq. \eqref{eq:cenersubs} can be calculated as
\begin{equation}
\bm{\sigma}_p^\ast = \left[\begin{array}{c}
                        \sigma_x \\
                        \sigma_y \\
                        \tau_{xy}
                      \end{array} \right]^\ast_p =
                      \left[\begin{array}{c}
                        - I_x (b_x) \\
                        - I_y (b_y)  \\
                        0
                      \end{array} \right].
\end{equation}
Above, $b_x$ and $b_y$ are the body force components. Furthermore,  $I_x(b_x)$ is an antiderivative of $b_x$ with respect to the coordinate $x$, i.e. it holds $\frac{\partial}{\partial x}I_x(b_x) = b_x$. Analogous definition applies to $I_y(b_y)$ with respect to the coordinate $y$. If the antiderivatives for $b_x$ or $b_y$ are not explicitly known, we first approximate the body force within the element by a constant field, then we take the corresponding antiderivatives.

Substituting Eq. \eqref{eq:stressModes} into Eq. \eqref{eq:cener} leads to
\begin{equation}
J^c(\bm{\beta}) = \frac{1}{2} \int_E (\bm{P}\bm{\beta} + \bm{\sigma}_p^\ast)^T \bm{C}^{-1} (\bm{P}\bm{\beta} + \bm{\sigma}_p^\ast) - \int_{\partial E} \bm{u}^T \bm{N}_E (\bm{P}\bm{\beta} + \bm{\sigma}_p^\ast).
\label{eq:cenersubs}
\end{equation}

Imposing null variations of the functional expressed in Eq. \eqref{eq:rcpfun}, the following is obtained
\begin{equation}
\bm{H} \bm{\beta} = \bm{g},
\label{eq:rcpfun}
\end{equation}
\noindent where
\begin{equation}
\bm{H} = \int_E \bm{P}^T \bm{C}^{-1} \bm{P}, \qquad \bm{g} = \int_{\partial E} \bm{P}^T \bm{N}_{E}^T \bm{u}  - \int_{E}\bm{P}^T\bm{C}^{-1}\bm{\sigma}_p^\ast. \label{eq:hg}
\end{equation}

\noindent In the definition of $\bm{g}$, see Eq. \eqref{eq:hg}, the first term can be integrated by parts and, exploiting that the stress modes are divergence-free, we get
\begin{equation}
\int_{\partial E} \bm{P}^T \bm{N}_{E}^T \bm{u}  = \int_{E} \bm{P}^T \bm{\varepsilon}(\bm{u}). \label{eq:g}
\end{equation}

In the finite element framework, the vector $\bm{g}$ can be computed using the left-hand side or the right-hand side of Eq. \eqref{eq:g}, leading to identical results. On the contrary, when virtual element schemes are considered, the left-hand side should be used, since the displacements are pointwise known only at the element boundary.
%
%

\begin{figure}[h!]
	\centering
	\includegraphics[width=0.80\textwidth,trim = 0mm 0mm 0mm 0mm, clip]{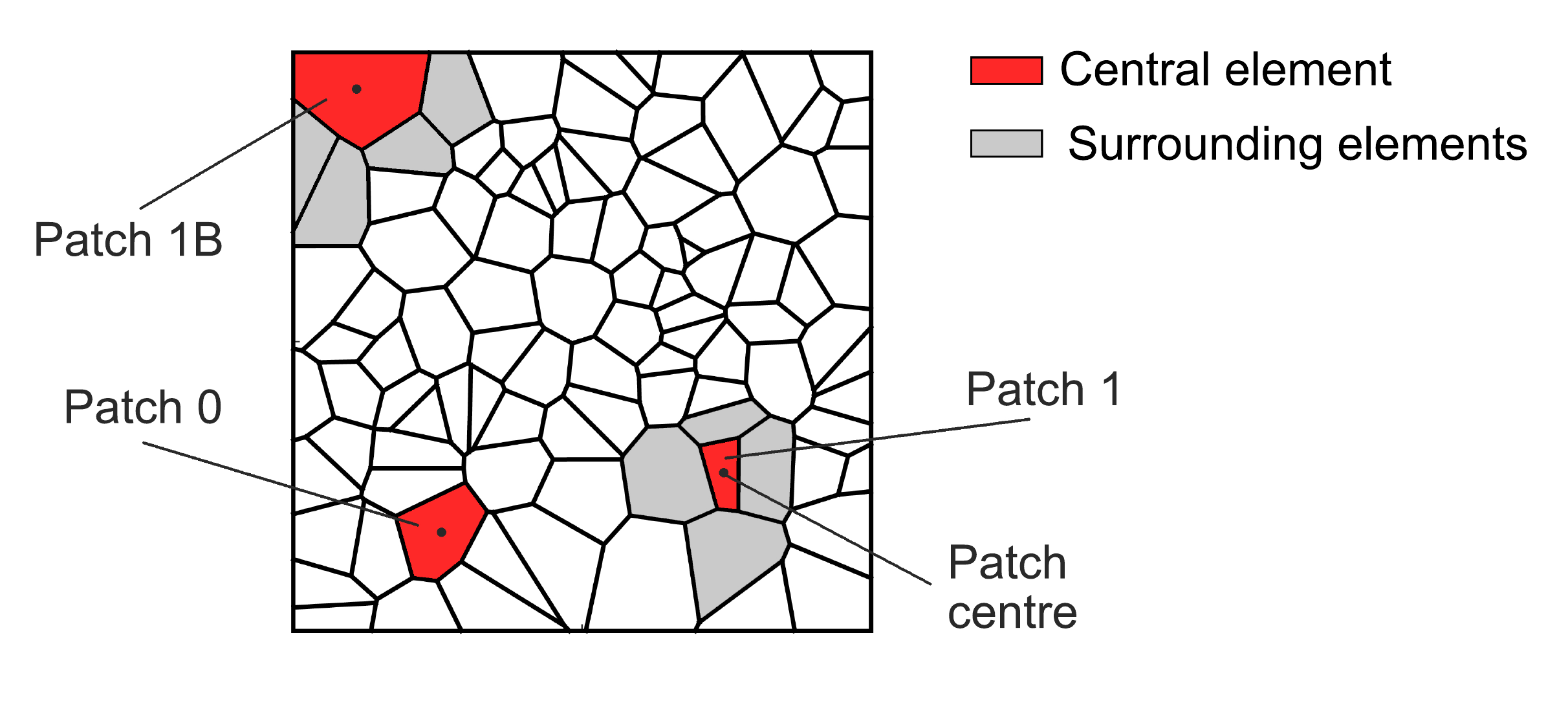}
	\caption{Patch types used in the present study.}
	\label{fig:patches}
\end{figure}

The adoption of a patch of elements instead of a single element does not introduce additional complexity: simply, the integrals appearing in Eq. \eqref{eq:hg} must be computed over the entire patch. Figure \ref{fig:patches} provides examples of the patches adopted in the following numerical tests. In particular, \textit{Patch 0} indicates the degenerate case of a single element, \textit{Patch 1} is a patch obtained by considering all surrounding elements sharing at least one node with a central element, while \textit{Patch 1B} represents a {\em boundary} patch of the same type as \textit{Patch 1}.

\section{Numerical results}
\label{sec:num}
In this section, numerical tests are used in order to validate the proposed stress recovery procedure. In particular, a square domain characterised by edge length equal to unity is considered. The body is assumed to be homogeneous and isotropic and its mechanical parameters are assigned in terms of the Lam\'e constants $\lambda = 1$ and $\mu = 1$ \cite{BeiraoLovaMora,Artioli2017a}. Plane strain regime is assumed.

Indicating as $\bm{\sigma}^{ex}$ the analytical stress field, the following error norm is used in order to evaluate the accuracy of the obtained results
\begin{equation}
    \label{eq:stress_err_norm}
    E_{\sigma} = \int_\Omega (\bm{\sigma}^{ex} - \bm{\sigma}^\ast)^T \bm{C}^{-1} (\bm{\sigma}^{ex} - \bm{\sigma}^\ast).
\end{equation}

\begin{figure}[h!]
	\centering
	\renewcommand{\thesubfigure}{}
	\subfigure[Tri (S)]{\includegraphics[width=0.24\textwidth,trim = 20mm 3mm 25mm 3mm, clip]{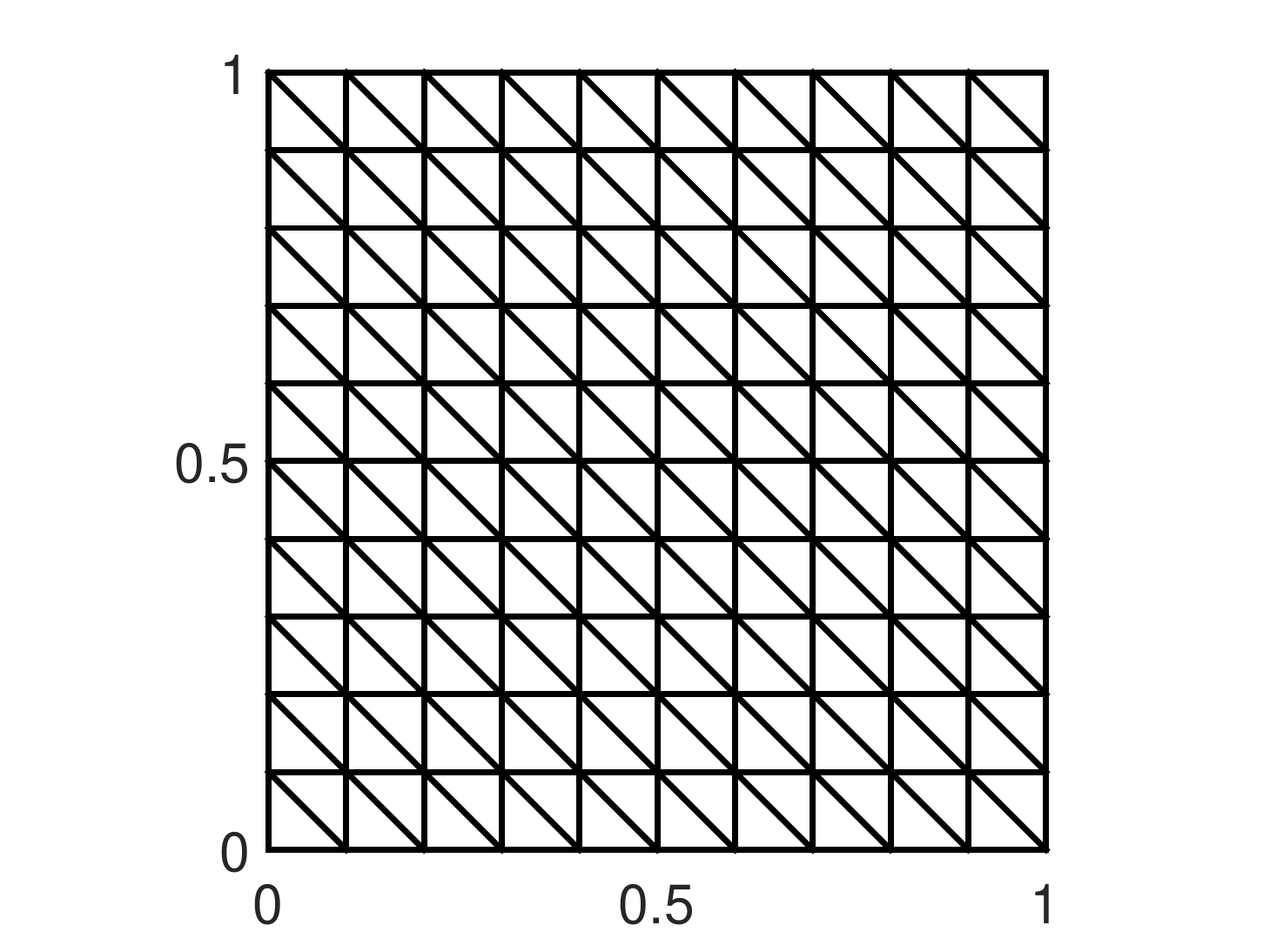}}
	\subfigure[Quad (S)]{\includegraphics[width=0.24\textwidth,trim = 20mm 3mm 25mm 3mm, clip]{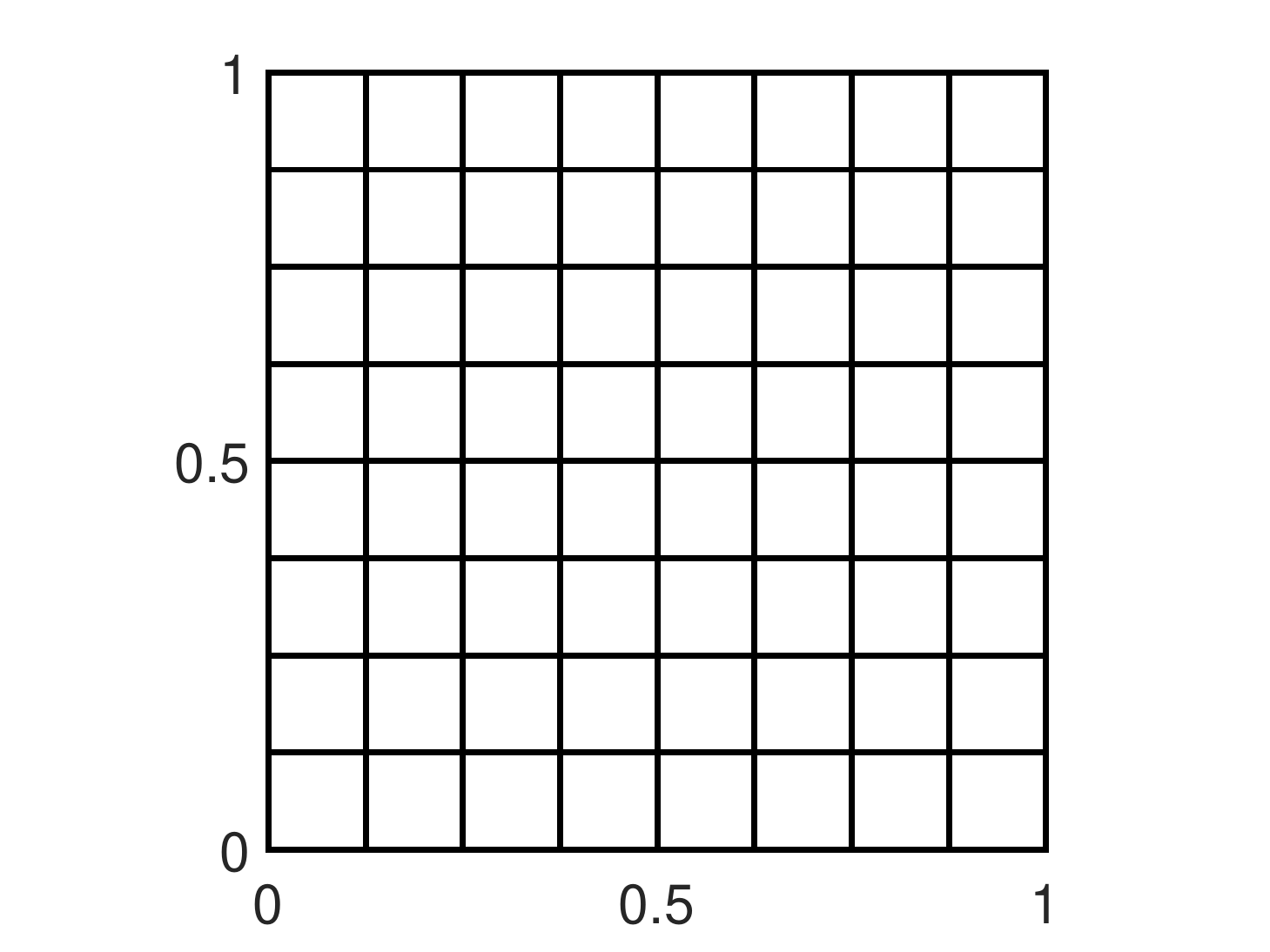}}
	\subfigure[Hex (S)]{\includegraphics[width=0.24\textwidth,trim = 20mm 3mm 25mm 3mm, clip]{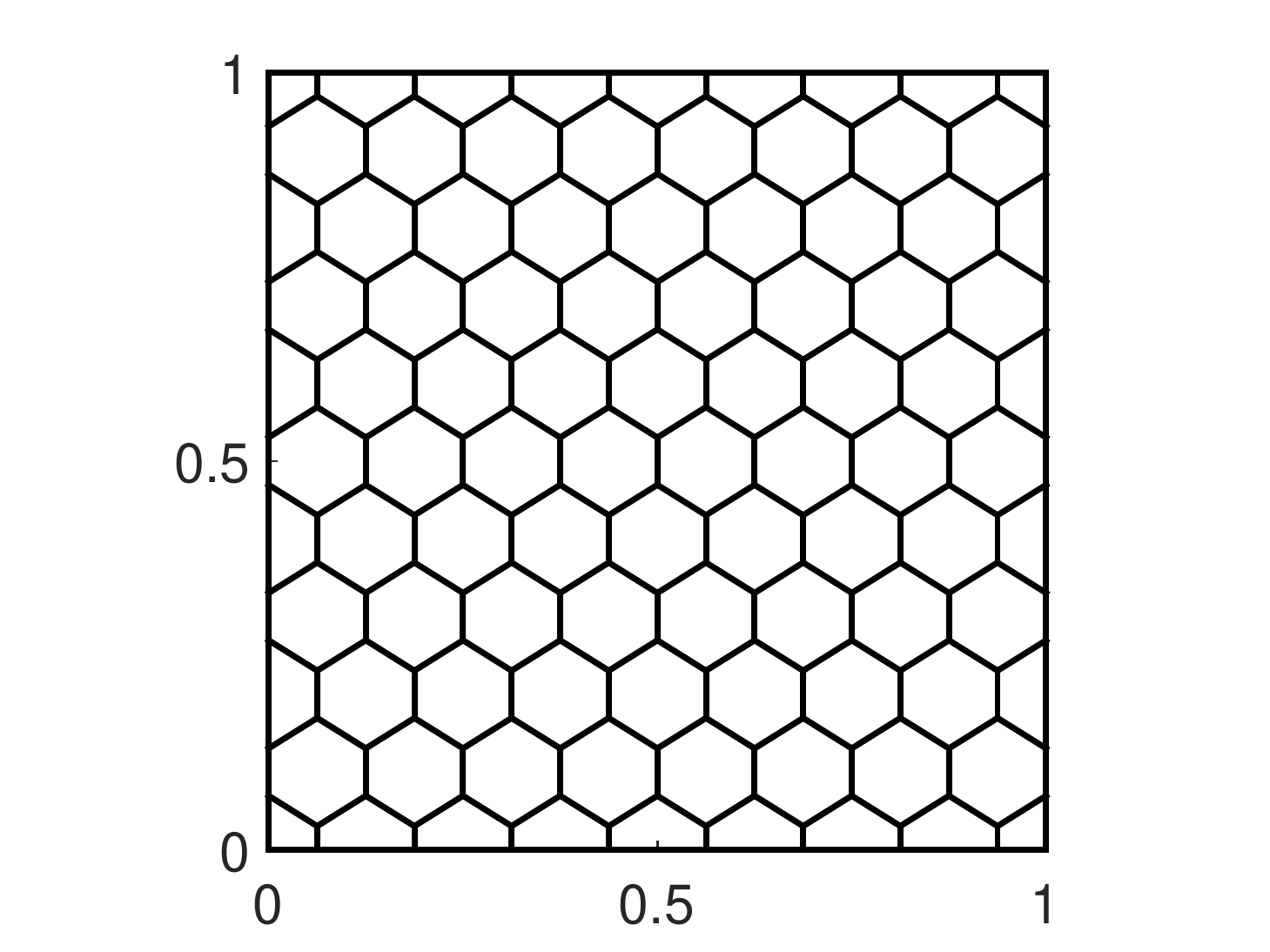}}
    \subfigure[Conc (S)]{\includegraphics[width=0.24\textwidth,trim = 20mm 3mm 25mm 3mm, clip]{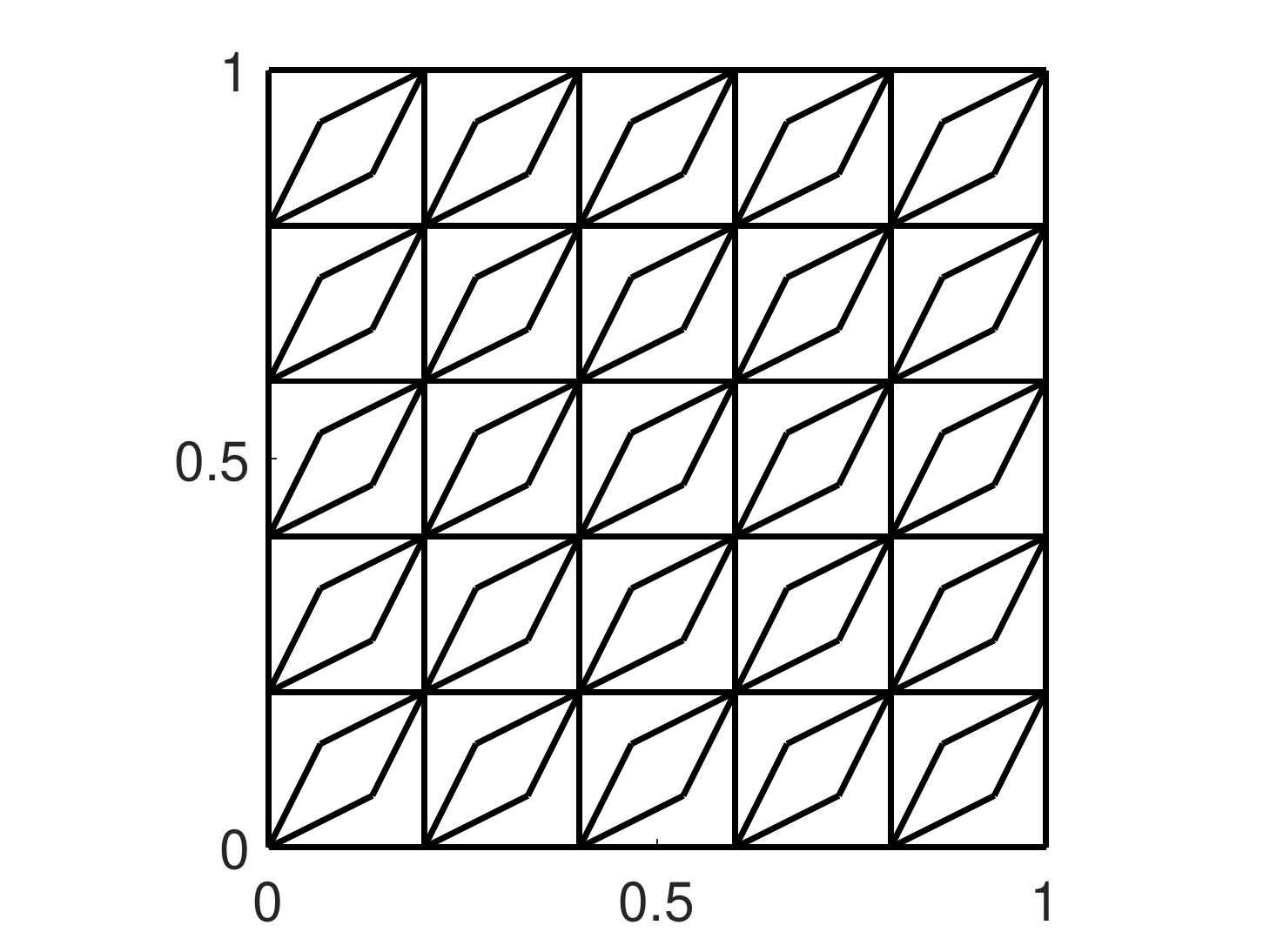}}
	\subfigure[Tri (U)]{\includegraphics[width=0.24\textwidth,trim = 20mm 3mm 25mm 3mm, clip]{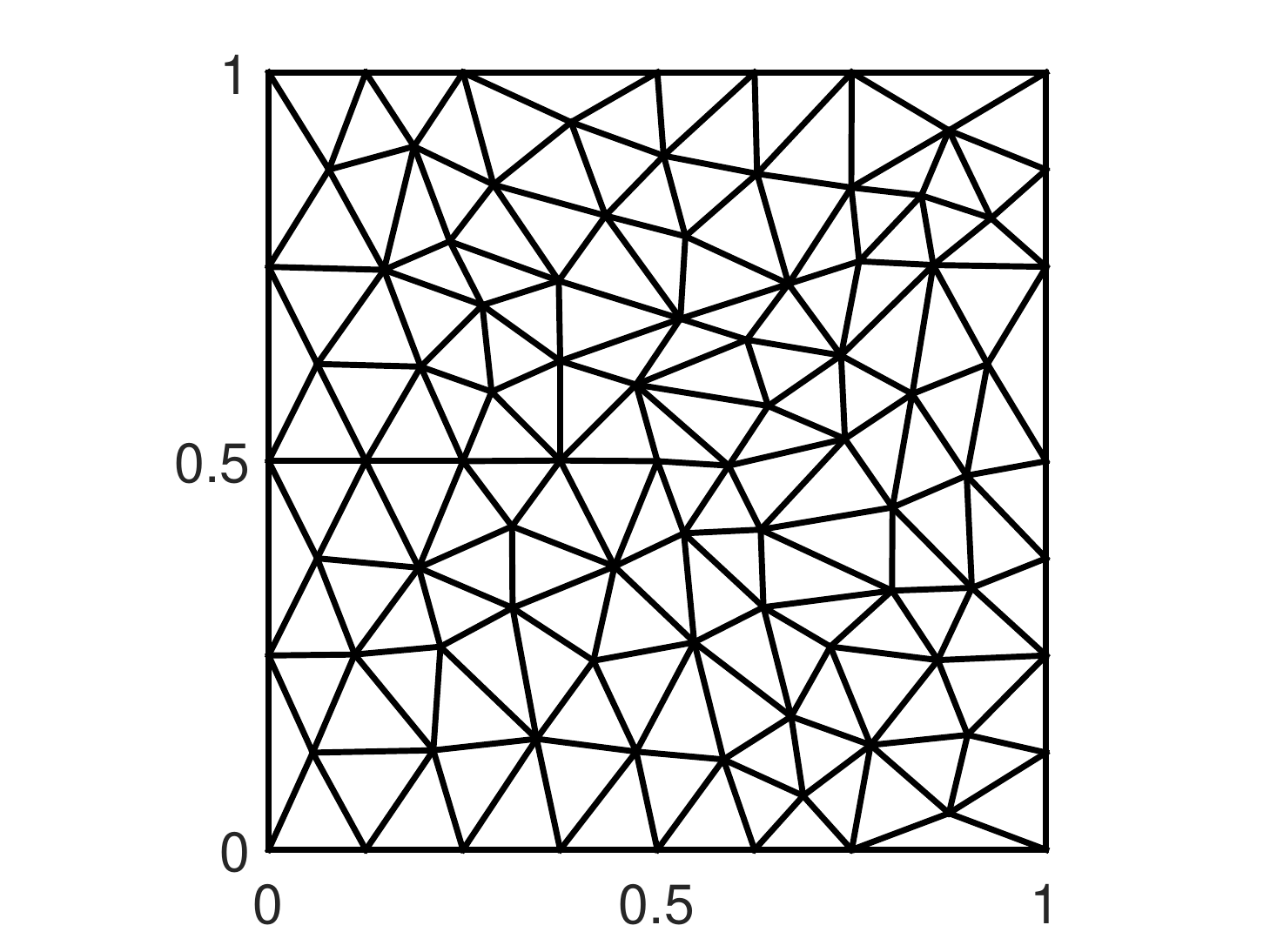}}
	\subfigure[Quad (U)]{\includegraphics[width=0.24\textwidth,trim = 20mm 3mm 25mm 3mm, clip]{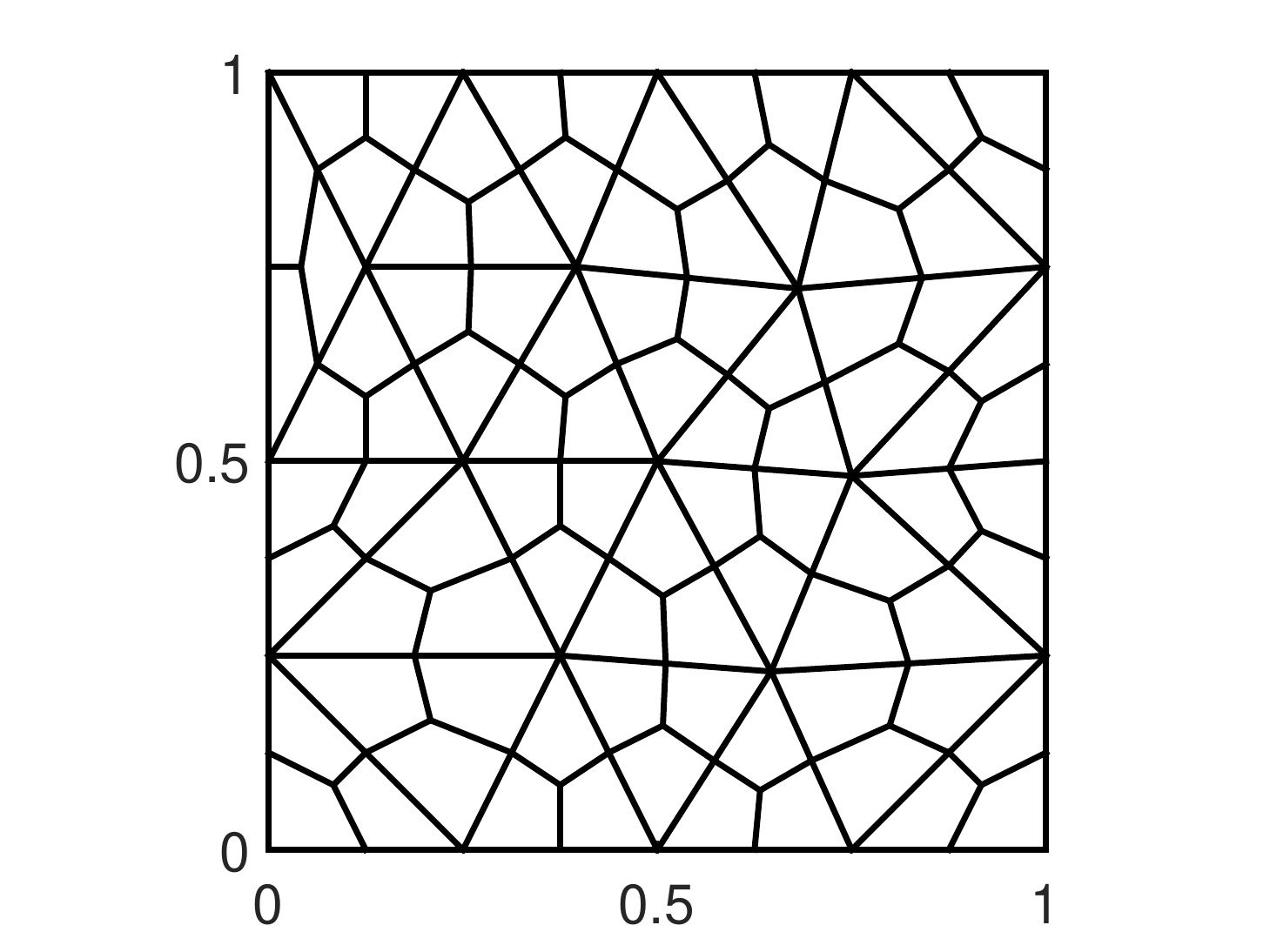}}
	\subfigure[Poly (U)]{\includegraphics[width=0.24\textwidth,trim = 20mm 3mm 25mm 3mm, clip]{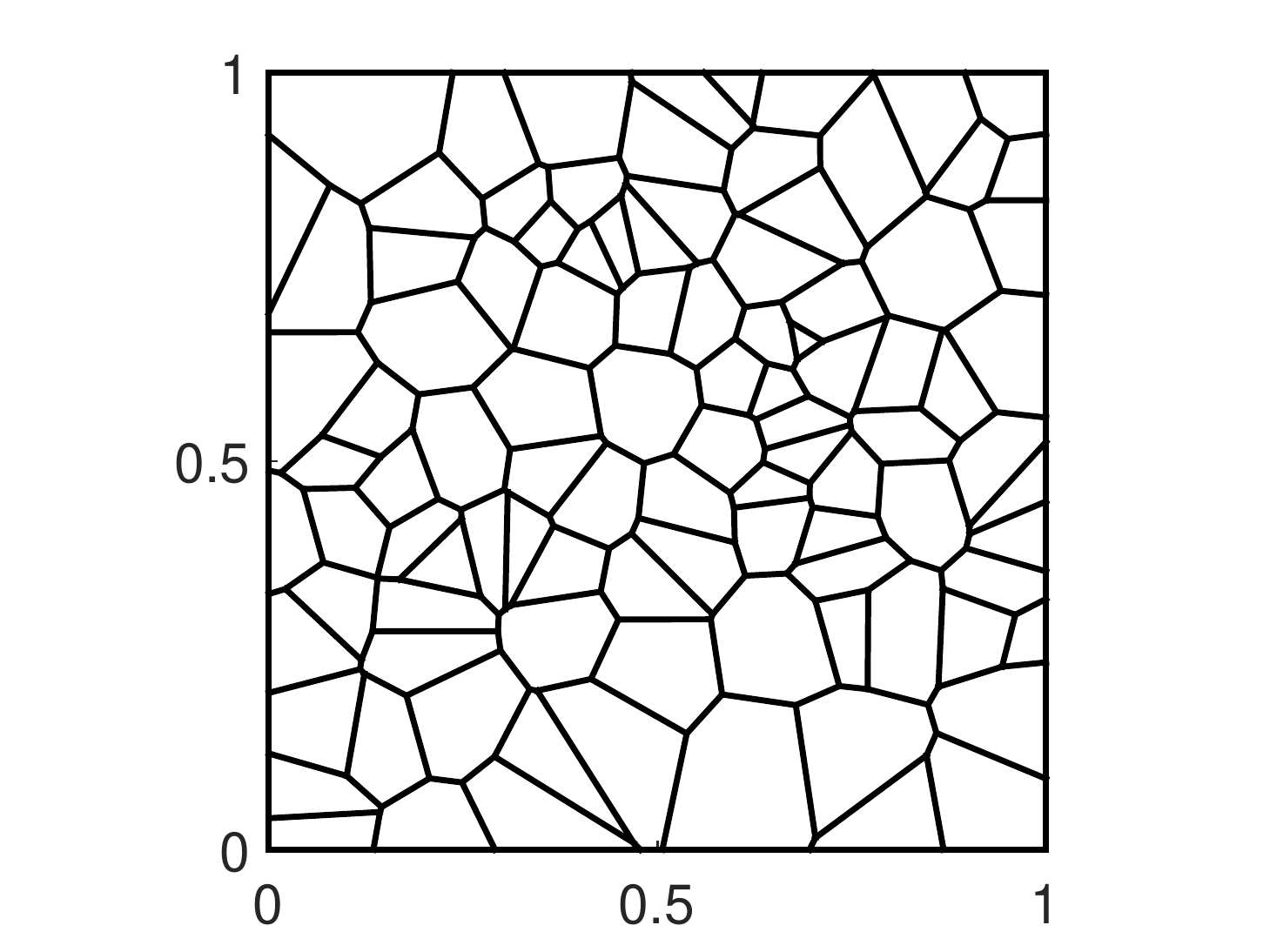}}
    \subfigure[Conc (U)]{\includegraphics[width=0.24\textwidth,trim = 20mm 3mm 25mm 3mm, clip]{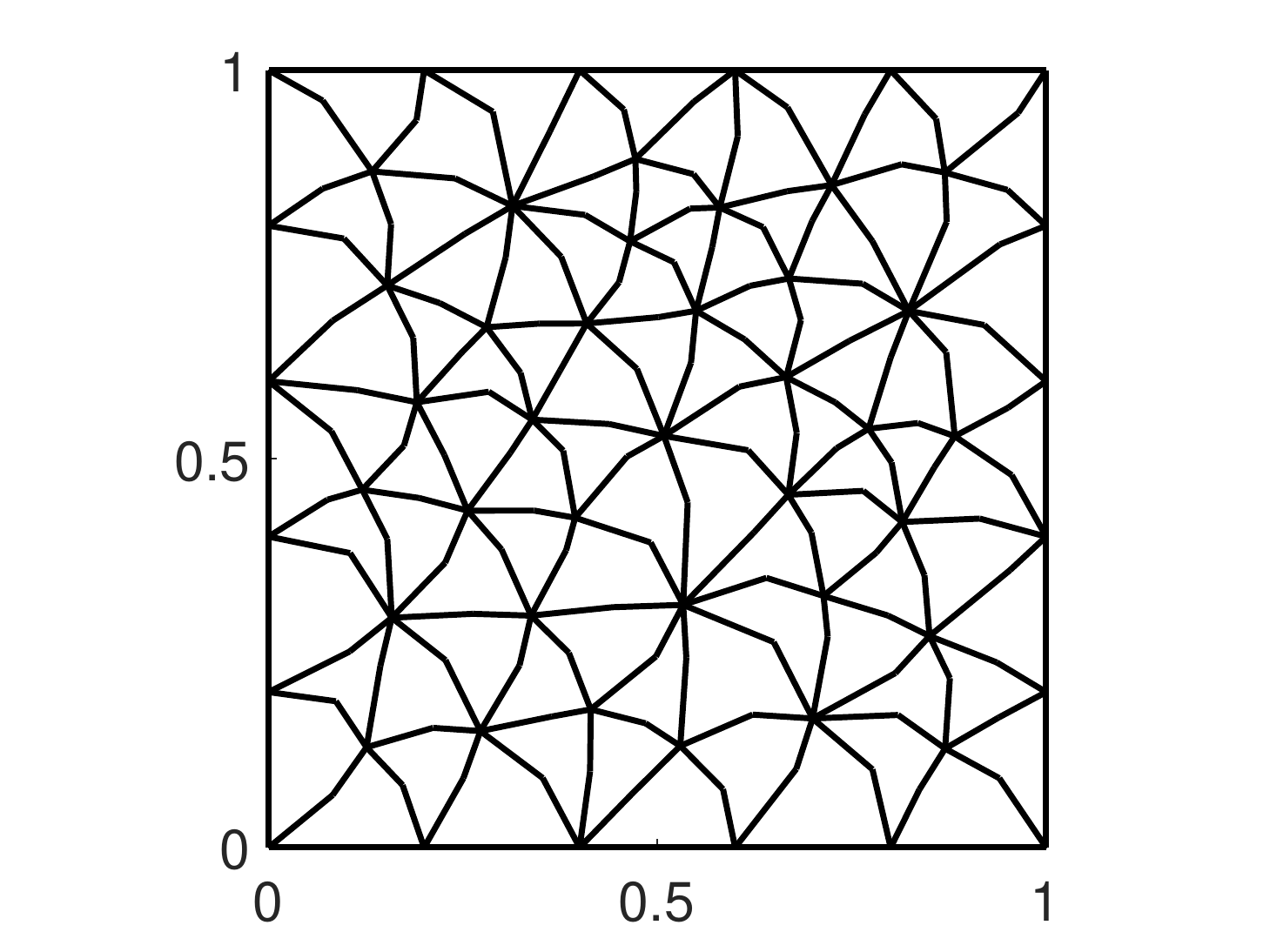}}
	\caption{Overview of adopted meshes for convergence assessment numerical tests \cite{Artioli2017HR}.}
	\label{fig:meshes}
\end{figure}

As it can be seen in Fig. \ref{fig:meshes}, eight meshes are considered. Four of them are structured and are indicated with the letter 'S' in the following. Those are composed of triangles, quadrilaterals, hexagons and a mixture of convex and concave quadrangles. The remaining four are unstructured meshes composed of triangles, quadrilaterals, random polygons and concave hexagons. Those meshes are indicated with the letter 'U'. The average edge length is indicated as $\bar{h}_e$.

Prescribed displacements are assumed and the corresponding body forces are calculated. Such body forces are then applied to the body, together with Dirichlet conditions over the entire boundary. Three tests are considered:
\begin{itemize}
\item  \textbf{Test a}: $u_x = x^3 - 3 x y^2, \; u_y = y^3 - 3 x^2 y$;
\item  \textbf{Test b}: $u_x = u_y = \sin(\pi x)  \sin(\pi y)$;
\item  \textbf{Test c}: $u_x = xy\sin(\pi x)\sin(\pi y)\; u_y = 0$.
\end{itemize}

In particular, \textit{Test a} is characterised by polynomial solution and engenders null body forces with inhomogeneous Dirichlet boundary conditions, while in \textit{Test b} the displacement solution is trigonometric and is characterised by homogeneous Dirichlet boundary conditions, see \cite{Artioli2017HR}. Analogously to \textit{Test b}, \textit{Test c} is characterised by homogeneous Dirichlet boundary conditions, and $u_x$ is a product of polynomial and trigonometric functions, while $u_y$ vanishes.

Three stress recovery procedures are compared and denoted as \textit{VEM}, \textit{RCP0} and \textit{RCP1}. In particular, \textit{VEM} corresponds to the standard procedure based on the adoption of the operator $\bm{\Pi}^m$, together with the use of the constitutive law. Instead, \textit{RCP0} and \textit{RCP1} denote the use of RCP as presented in Eq. \eqref{eq:hg}, for \textit{Patch 0} and \textit{Patch 1}, respectively.

\begin{figure}[h!]
	\centering
	\subfigure[]{\includegraphics[width=0.48\textwidth,trim = 0mm 0mm 0mm 0mm, clip]{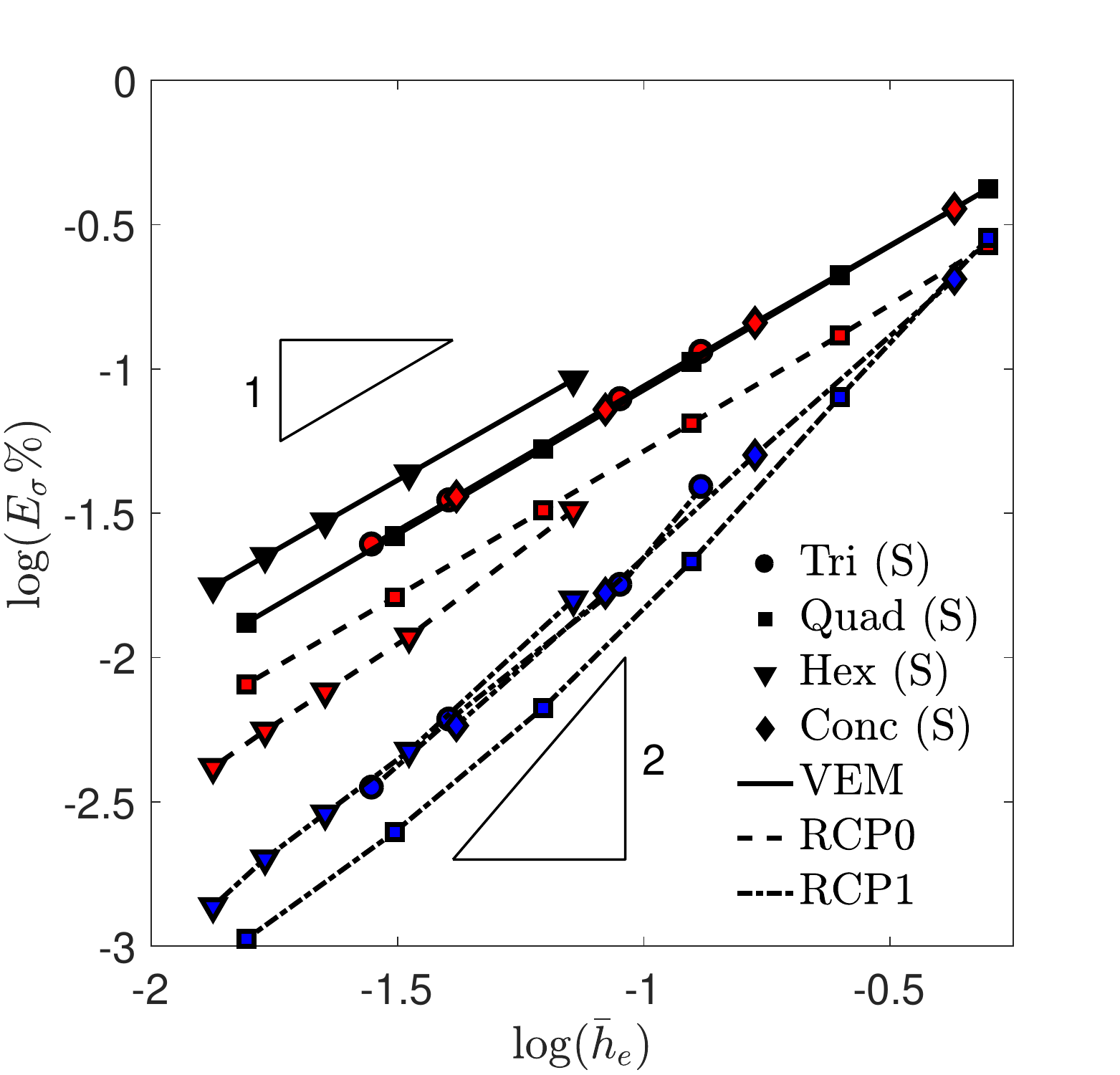}}
	\subfigure[]{\includegraphics[width=0.48\textwidth,trim = 0mm 0mm 0mm 0mm, clip]{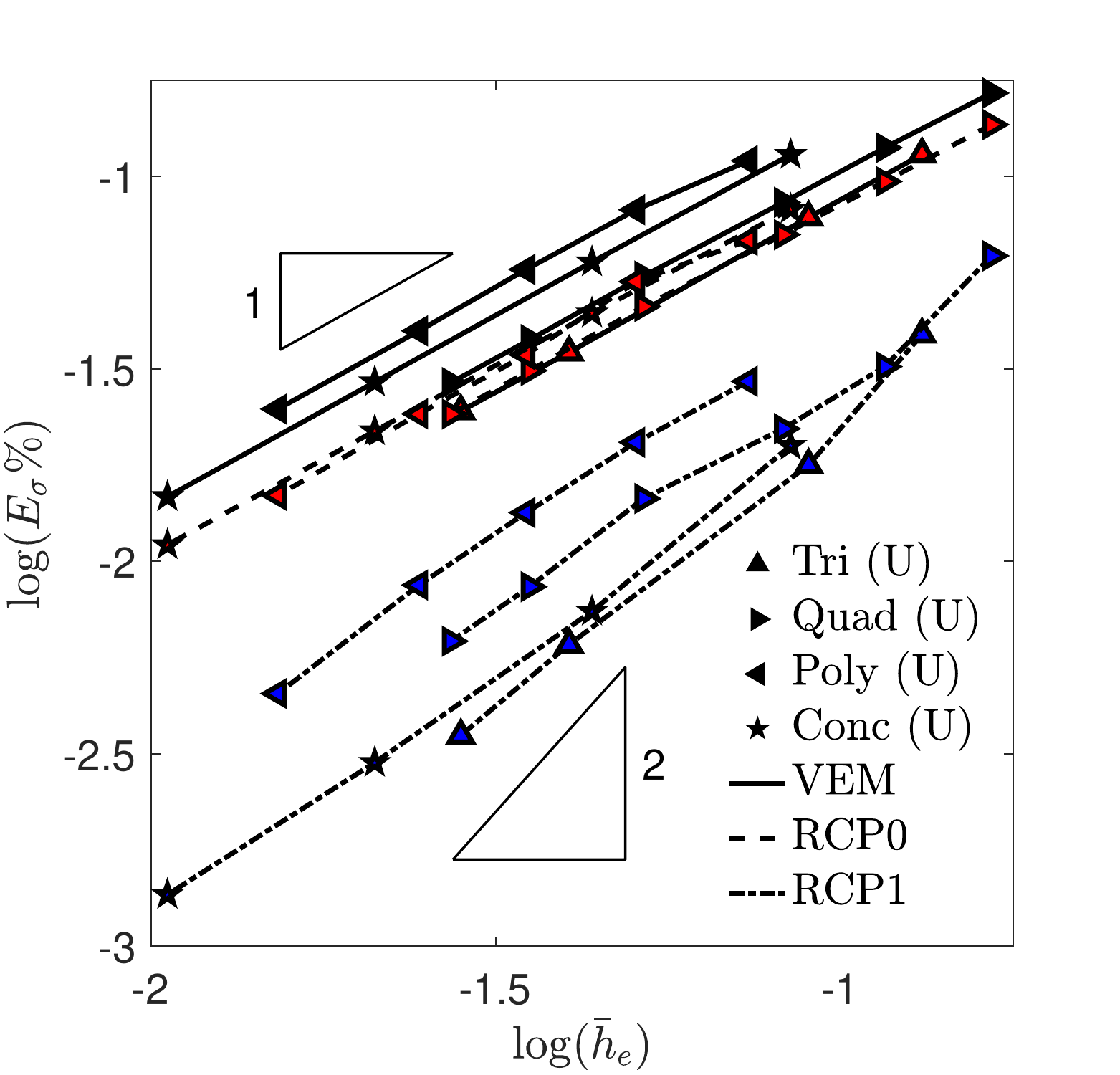}}\\
	\caption{$\bar{h}_e-$convergence results for Test $a$: (a) structured and (b) unstructured meshes. Black, red and blue markers indicate \textit{VEM}, \textit{RCP0} and \textit{RCP1}, respectively.}
	\label{fig:resuTestA}
\end{figure}

\begin{figure}[h!]
	\centering
	\subfigure[]{\includegraphics[width=0.48\textwidth,trim = 0mm 0mm 0mm 0mm, clip]{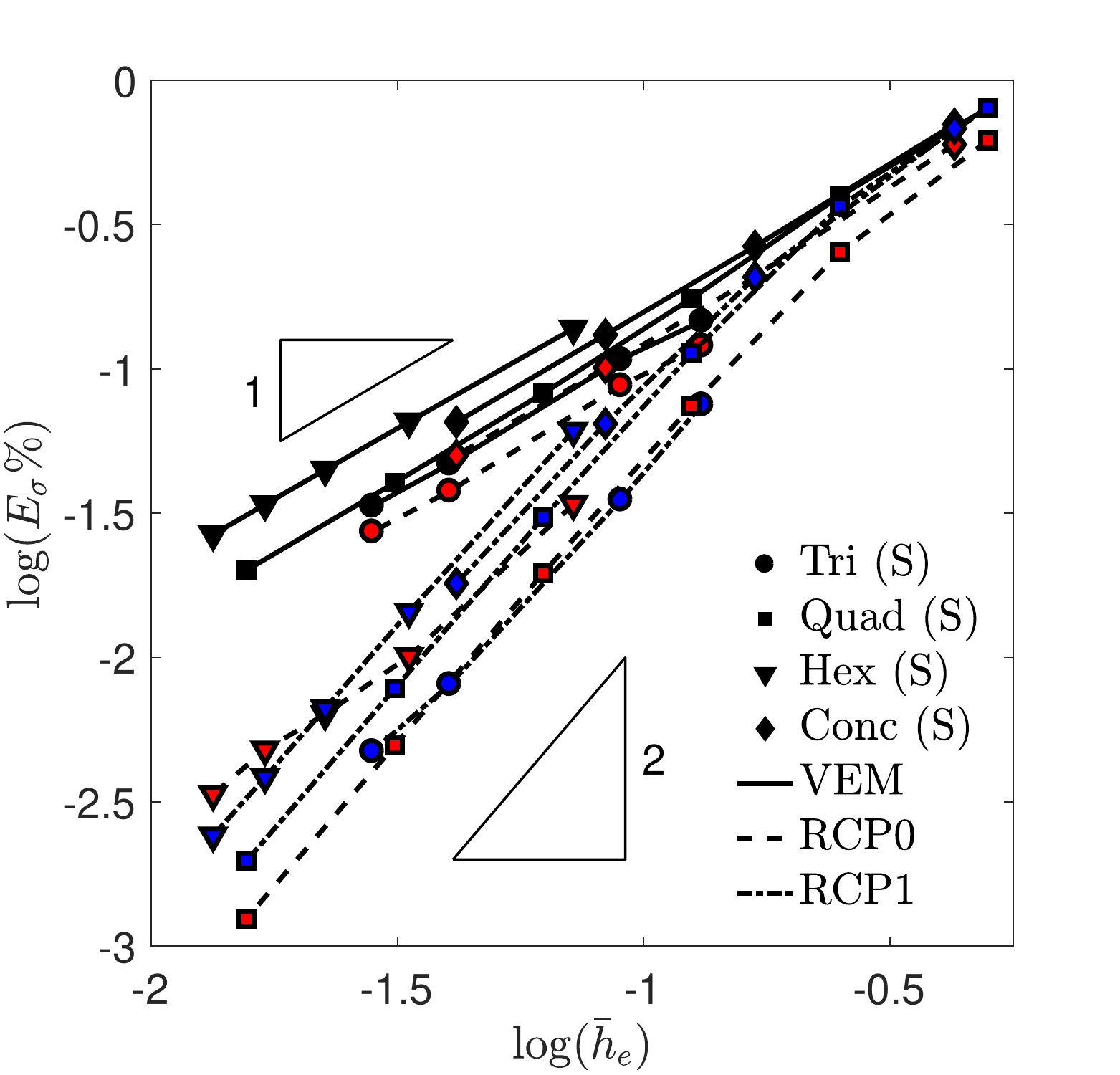}}
	\subfigure[]{\includegraphics[width=0.48\textwidth,trim = 0mm 0mm 0mm 0mm, clip]{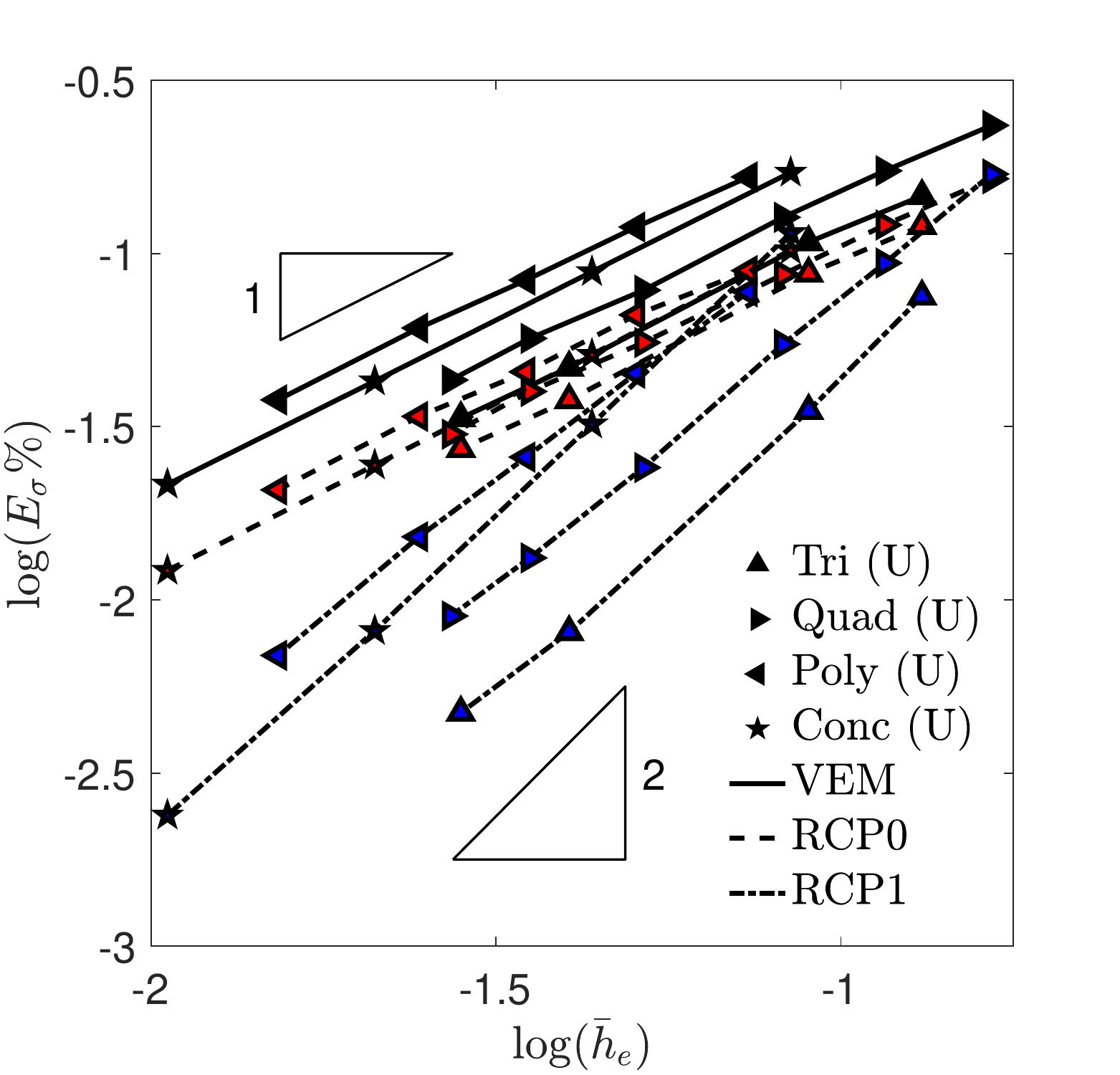}}\\
	\caption{$\bar{h}_e-$convergence results for Test $b$: (a) structured and (b) unstructured meshes. Black, red and blue markers indicate \textit{VEM}, \textit{RCP0} and \textit{RCP1}, respectively.}
	\label{fig:resuTestB}
\end{figure}

\begin{figure}[h!]
	\centering
	\subfigure[]{\includegraphics[width=0.48\textwidth,trim = 0mm 0mm 0mm 0mm, clip]{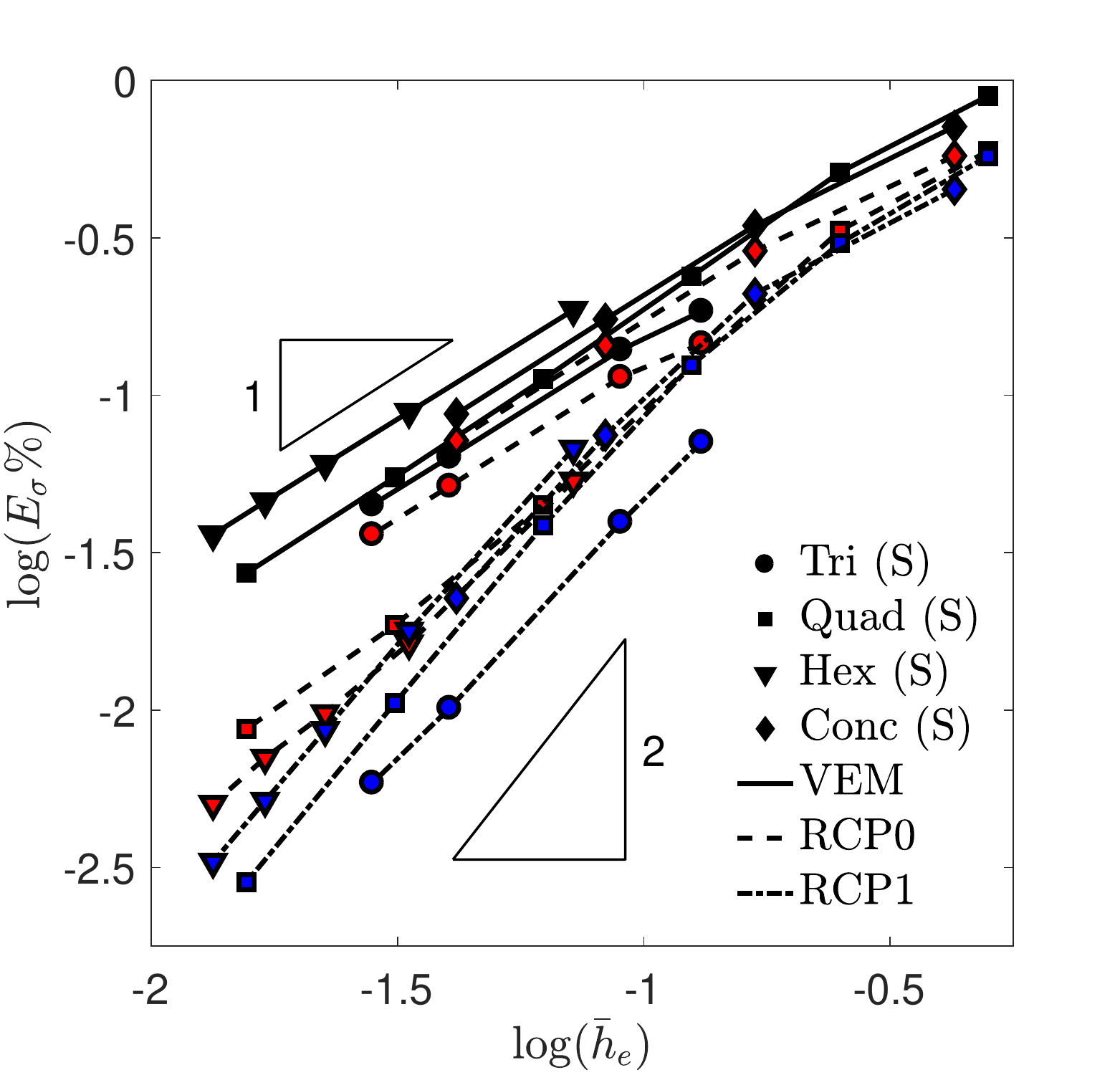}}
	\subfigure[]{\includegraphics[width=0.48\textwidth,trim = 0mm 0mm 0mm 0mm, clip]{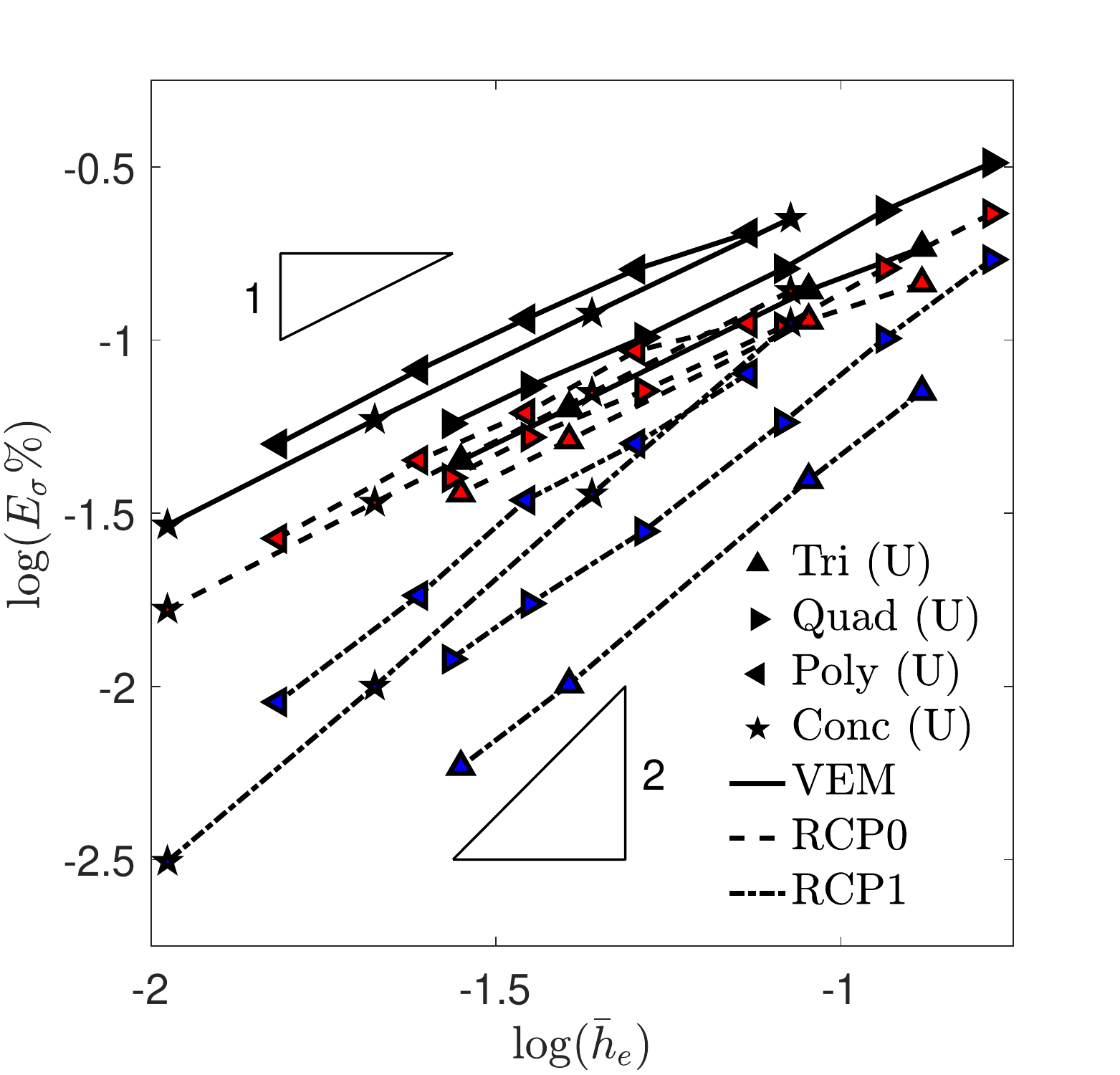}}\\
	\caption{$\bar{h}_e-$convergence results for Test $c$: (a) structured and (b) unstructured meshes. Black, red and blue markers indicate \textit{VEM}, \textit{RCP0} and \textit{RCP1}, respectively.}
	\label{fig:resuTestC}
\end{figure}

Figures \ref{fig:resuTestA}, \ref{fig:resuTestB} and \ref{fig:resuTestC} report results obtained for \textit{Test a}, \textit{Test b} and \textit{Test c} for all the aforementioned meshes, respectively. It can be clearly seen that RCP always has a favourable effect, leading to more accurate results if compared to the use of the $\bm{\Pi}^m$ operator and, in some cases, to a substantial increase in the convergence rate.
Regarding \textit{RCP0}, we notice that the improvement can be clearly appreciated when the Hex (S) mesh is considered while it tends to disappear for triangular elements.
Furthermore, considering now results obtained using \textit{RCP1}, it can be seen that the use of element patches has a significant beneficial effect, in agreement with what has been observed in the context of displacement-based finite elements \cite{Ubertini2004,Ubertini2006}.

\begin{figure}[h!]
	\centering
	\subfigure[]{\includegraphics[height=0.38\textwidth,trim = 0mm 20mm 10mm 20mm, clip]{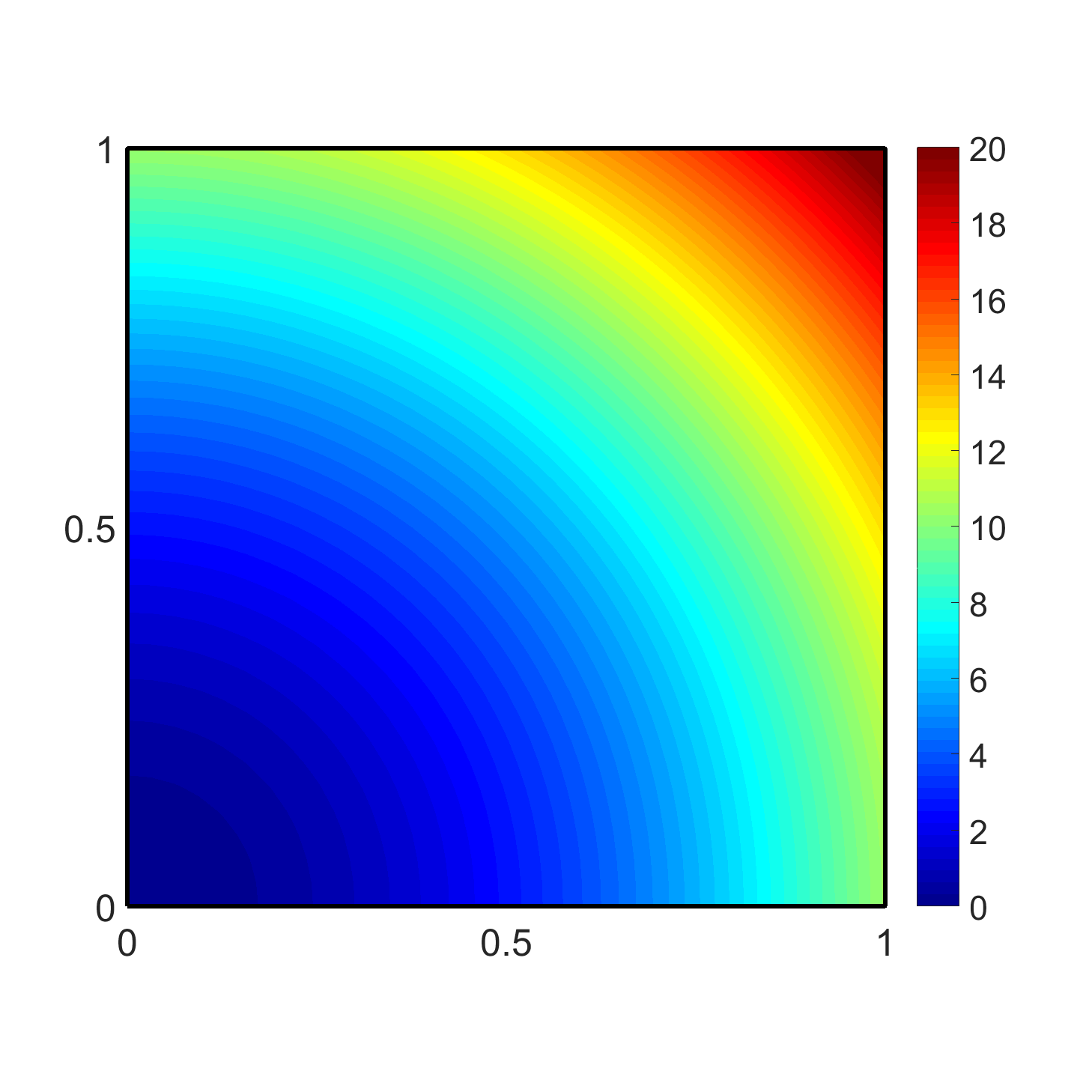}} \qquad
	\subfigure[]{\includegraphics[height=0.38\textwidth,trim = 0mm 20mm 10mm 20mm, clip]{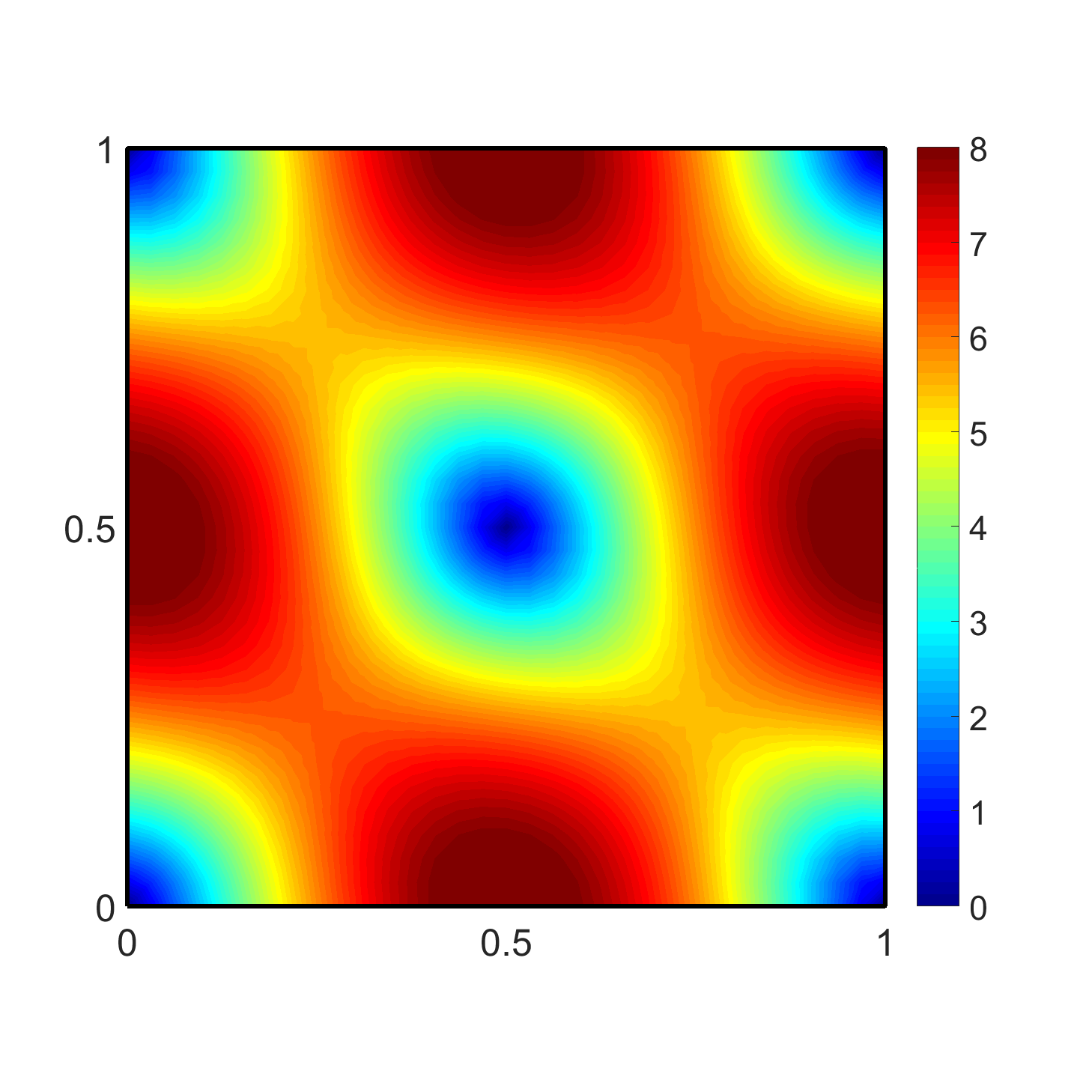}}
	\caption{Exact von Mises equivalent stress distributions: (a) \textit{Test a} and (b) \textit{Test b}.}
\label{fig:misesAnaSol}
\end{figure}

\begin{figure}[h!]
	\centering
	\subfigure[]{\includegraphics[height=0.3\textwidth,trim = 0mm 20mm 35mm 20mm, clip]{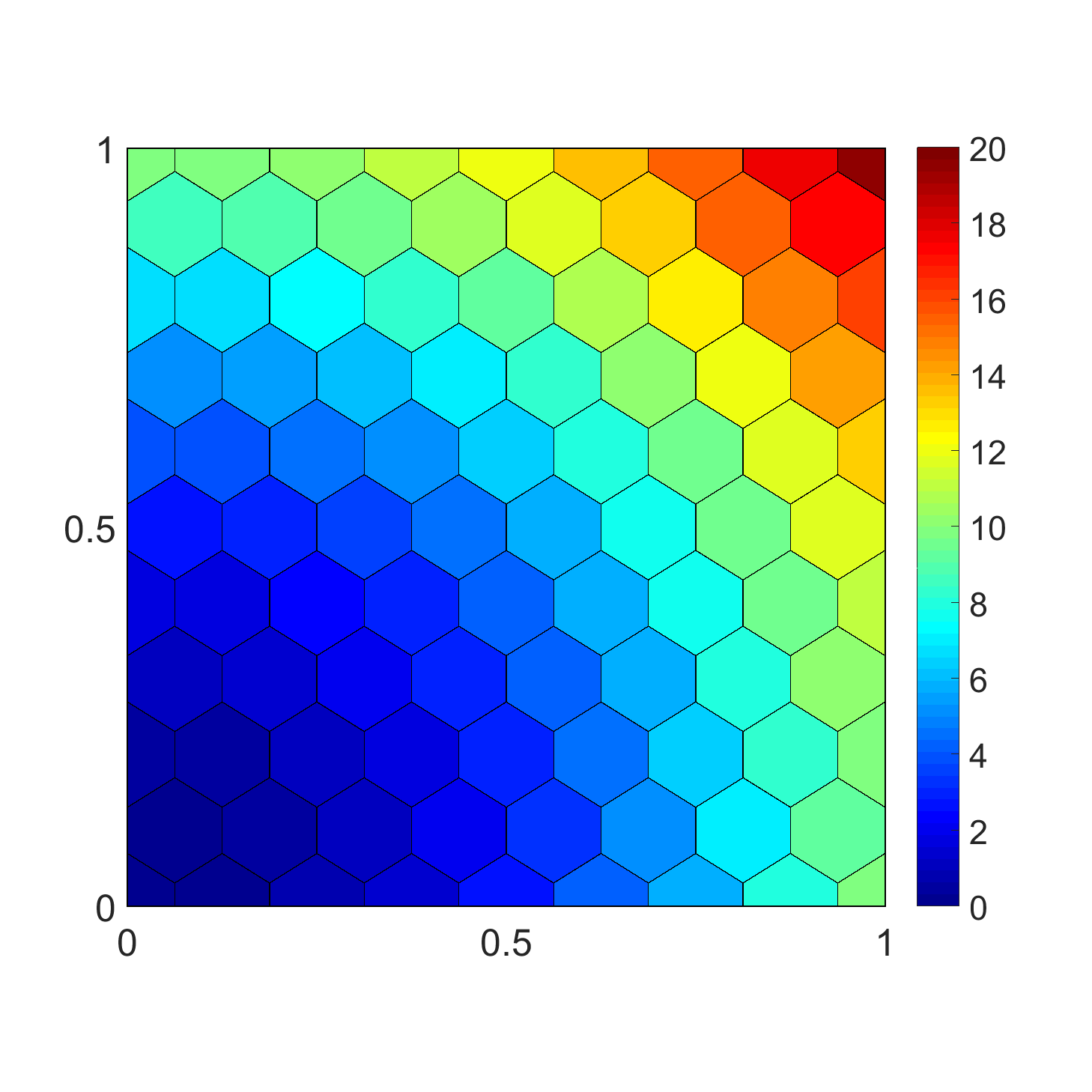}}
	\subfigure[]{\includegraphics[height=0.3\textwidth,trim = 0mm 20mm 35mm 20mm, clip]{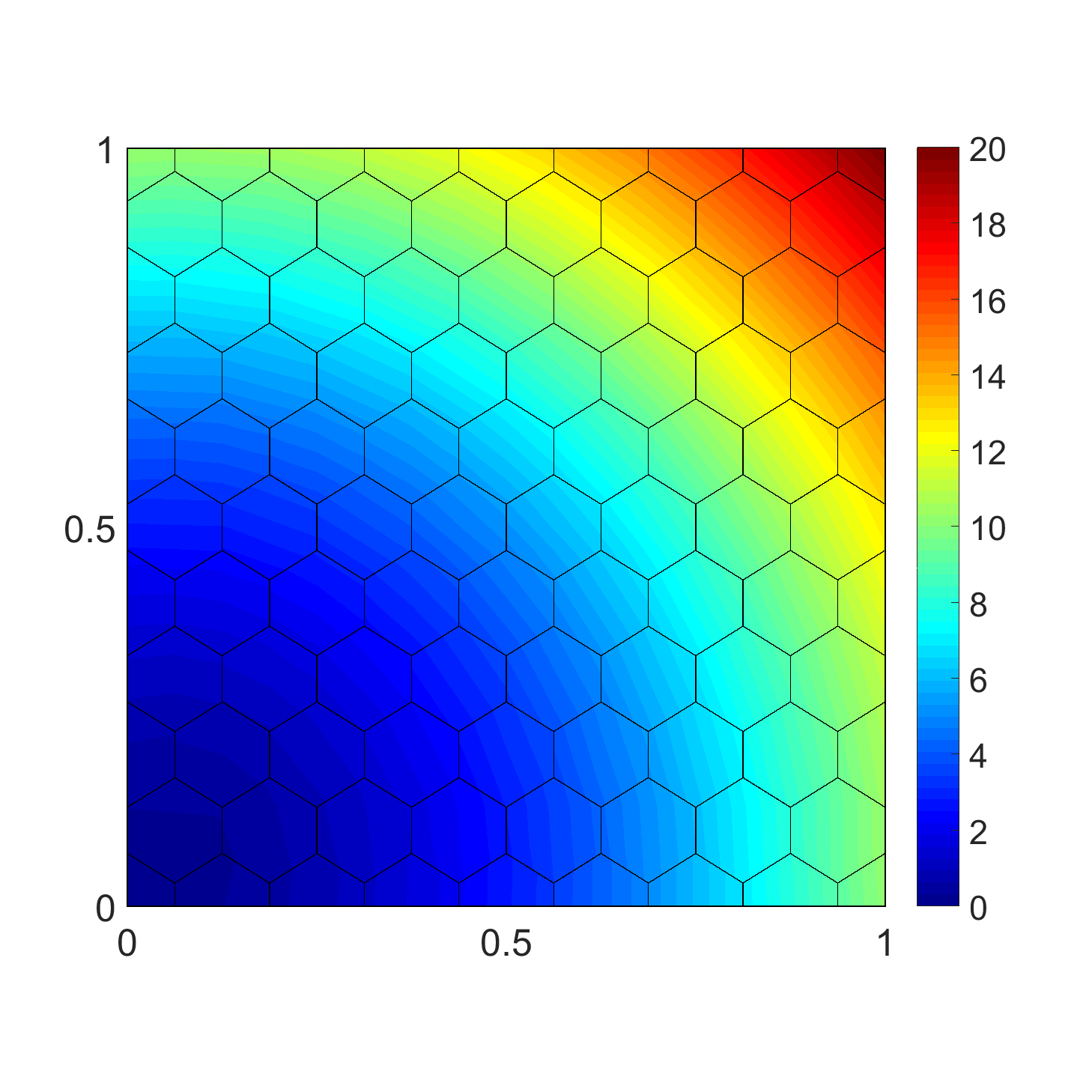}}
	\subfigure[]{\includegraphics[height=0.3\textwidth,trim = 0mm 20mm 10mm 20mm, clip]{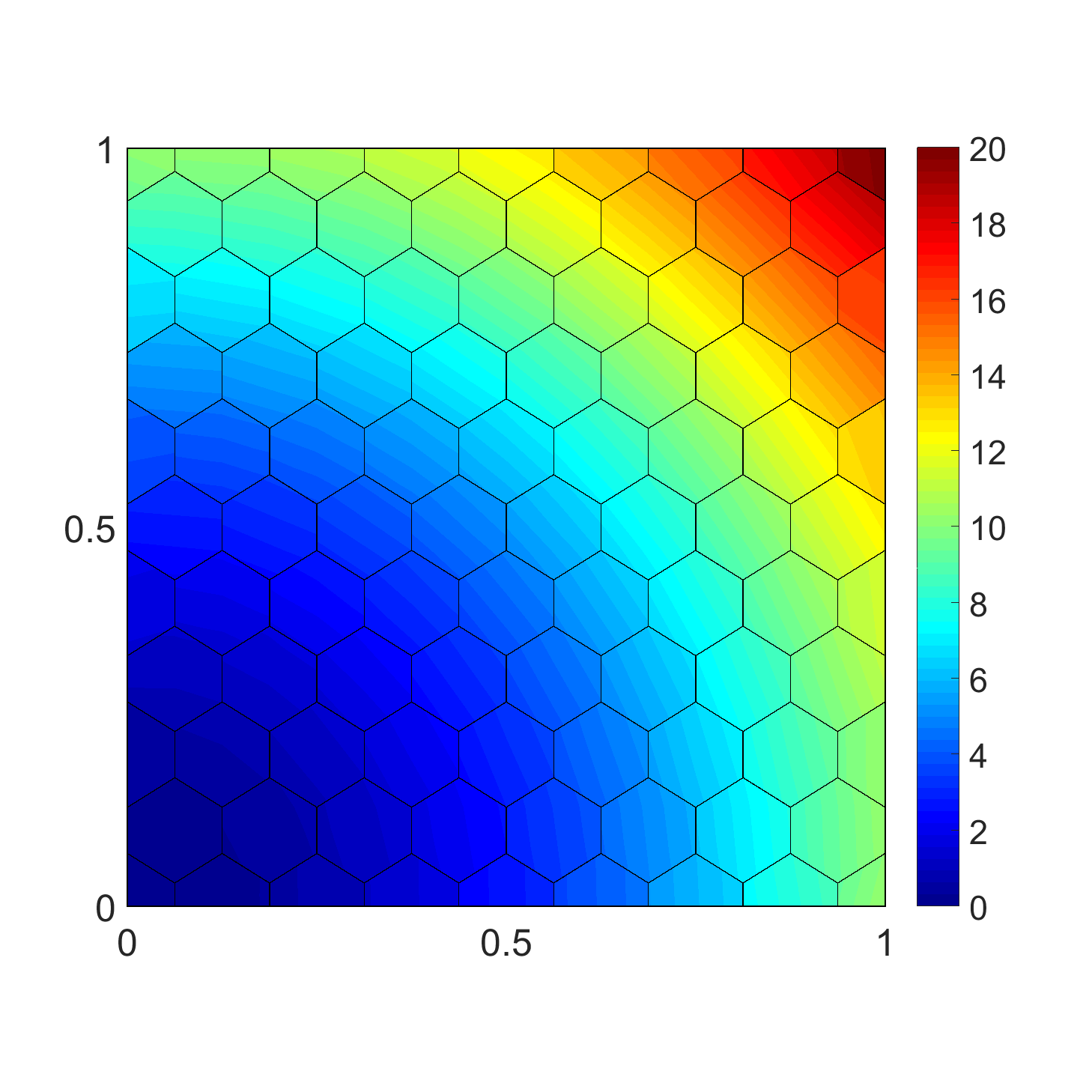}}\\
	\subfigure[]{\includegraphics[height=0.3\textwidth,trim = 0mm 20mm 35mm 20mm, clip]{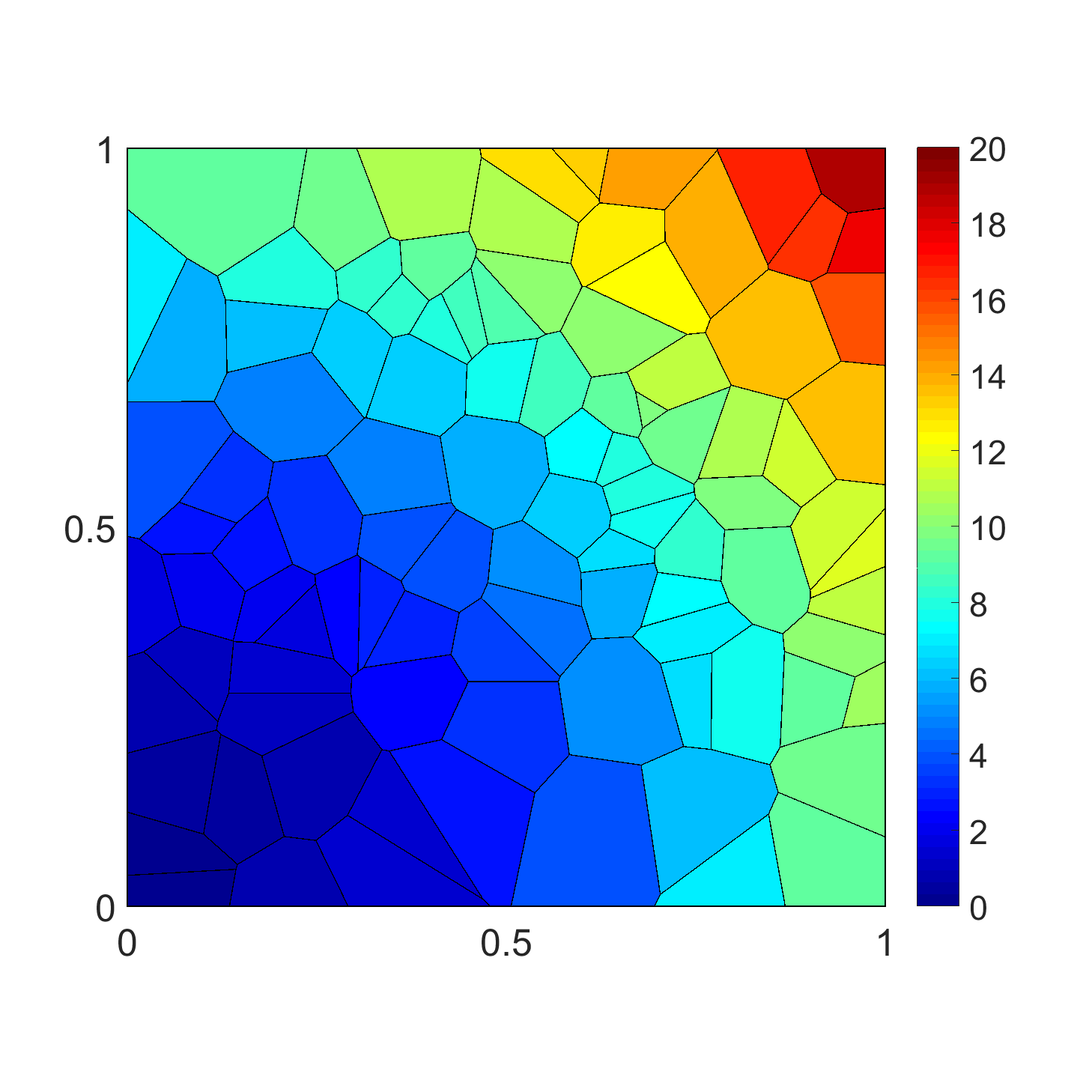}}
	\subfigure[]{\includegraphics[height=0.3\textwidth,trim = 0mm 20mm 35mm 20mm, clip]{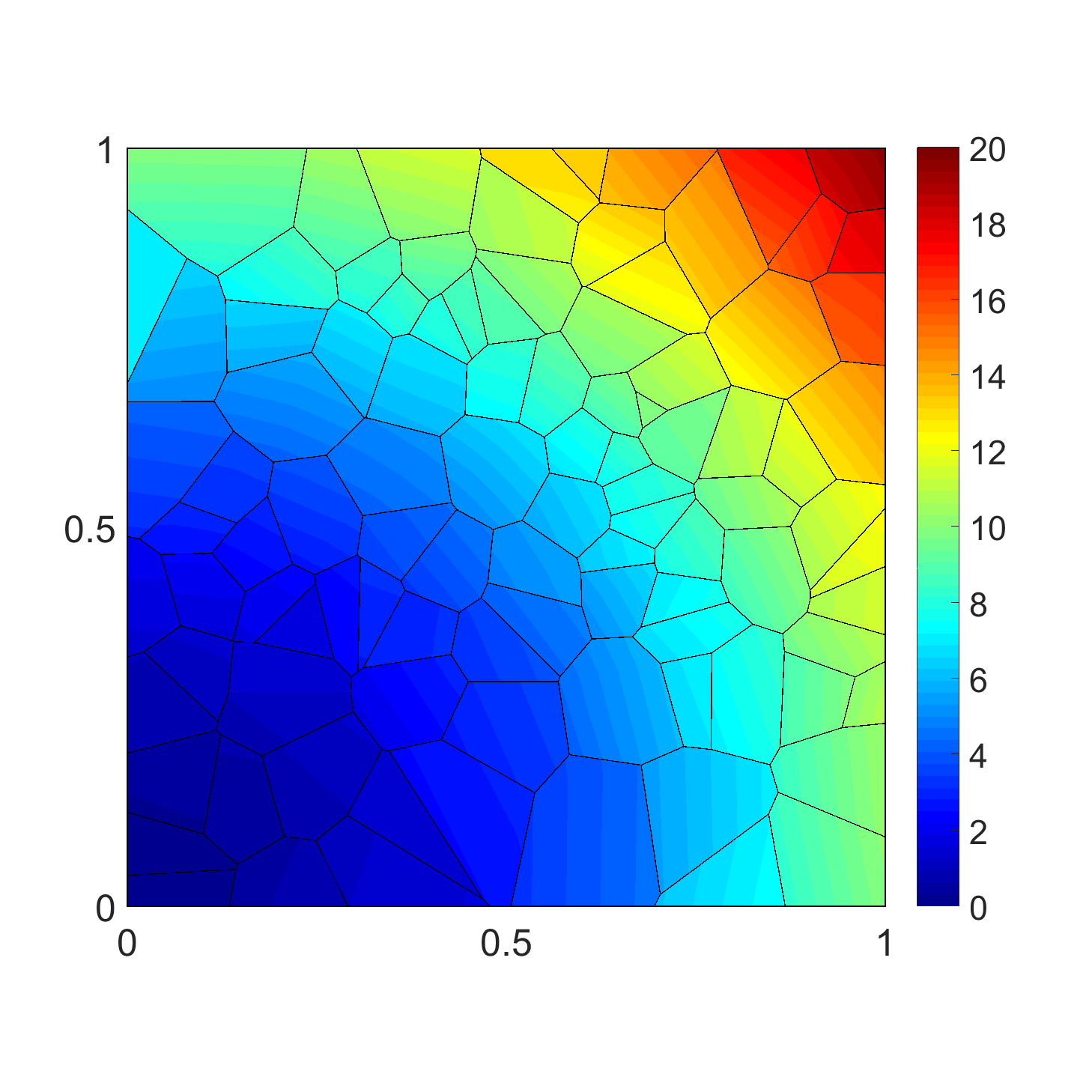}}
	\subfigure[]{\includegraphics[height=0.3\textwidth,trim = 0mm 20mm 10mm 20mm, clip]{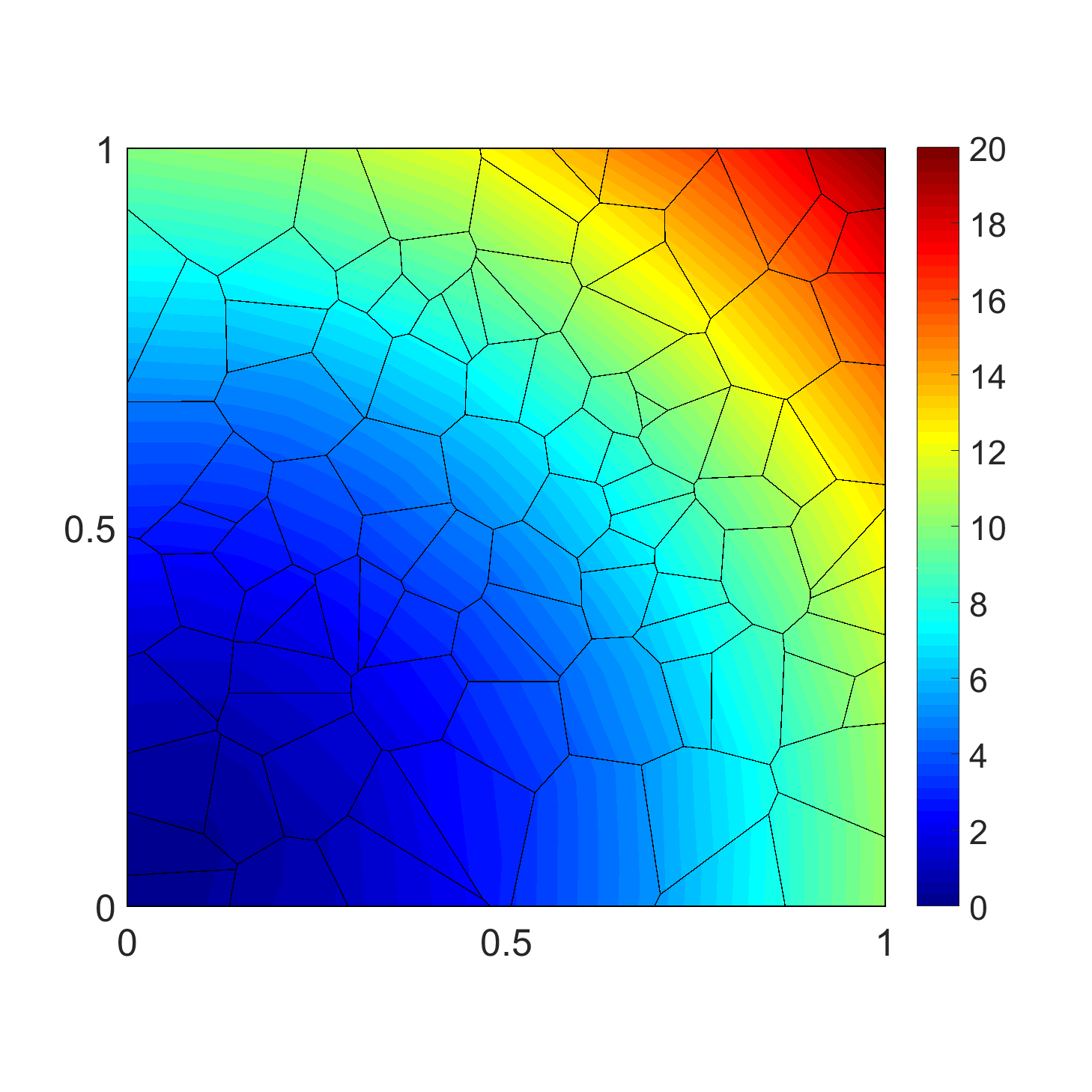}}\\
	\caption{Von Mises equivalent stress distributions for \textit{Test a}: (a) Hex (S) \textit{VEM}, (b) Hex (S) \textit{RCP0}, (c) Hex (S) \textit{RCP1}, (d) Poly (U) \textit{VEM}, (e) Poly (U) \textit{RCP0} and (f) Poly (U) \textit{RCP1}.}
	\label{fig:misesHexaTestA}
\end{figure}

\begin{figure}[h!]
	\centering
	\subfigure[]{\includegraphics[height=0.3\textwidth,trim = 0mm 20mm 35mm 20mm, clip]{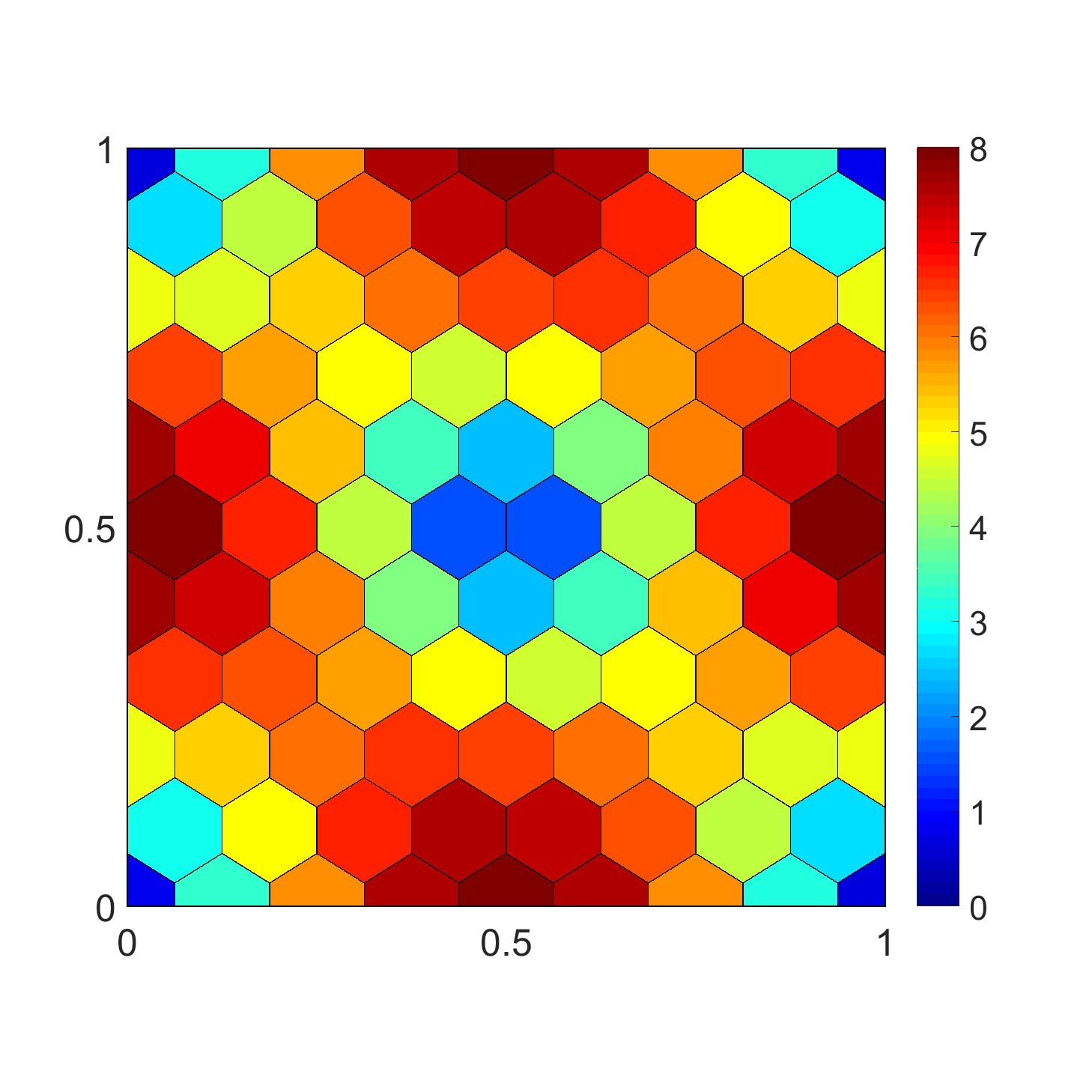}}
	\subfigure[]{\includegraphics[height=0.3\textwidth,trim = 0mm 20mm 35mm 20mm, clip]{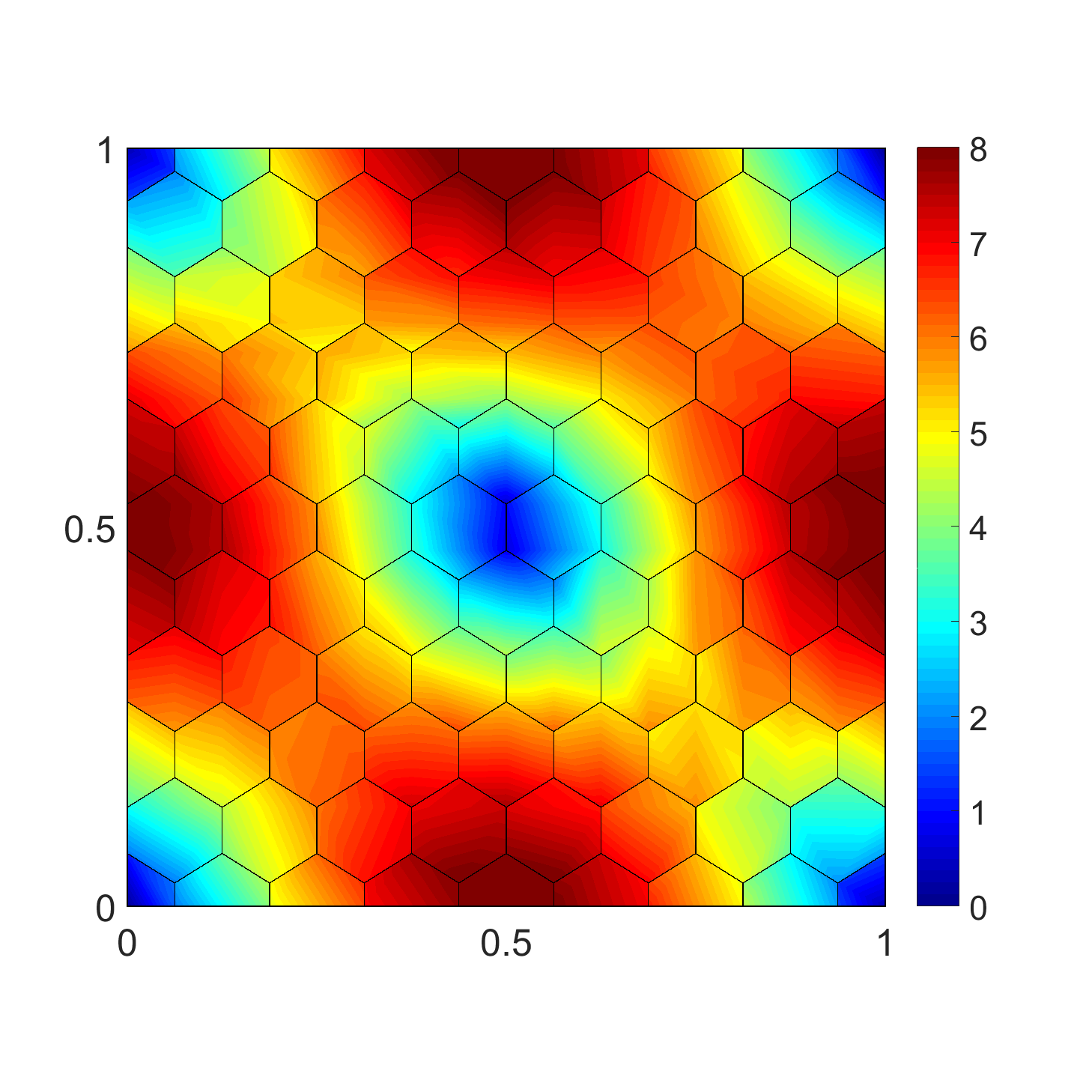}}
	\subfigure[]{\includegraphics[height=0.3\textwidth,trim = 0mm 20mm 10mm 20mm, clip]{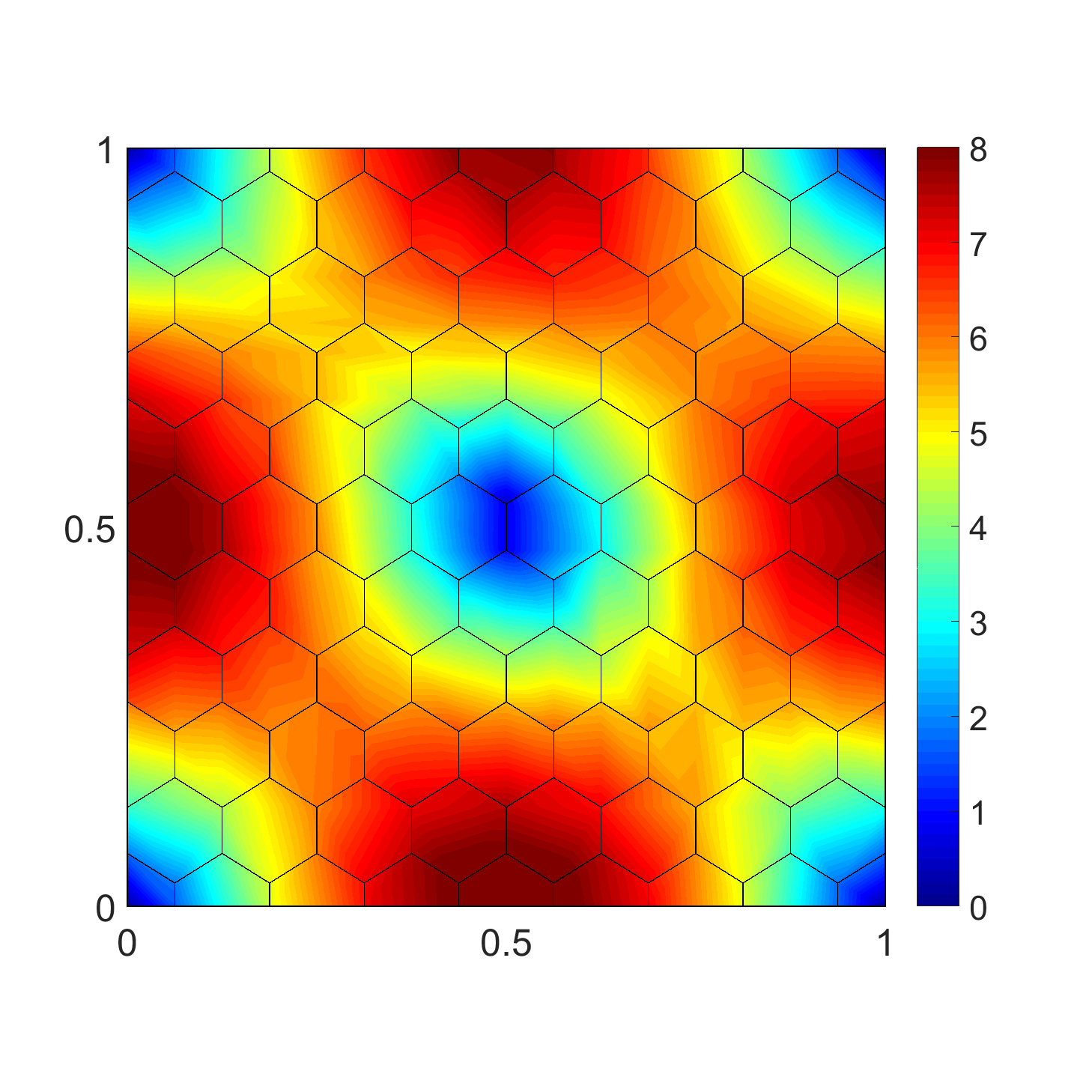}}\\
	\subfigure[]{\includegraphics[height=0.3\textwidth,trim = 0mm 20mm 35mm 20mm, clip]{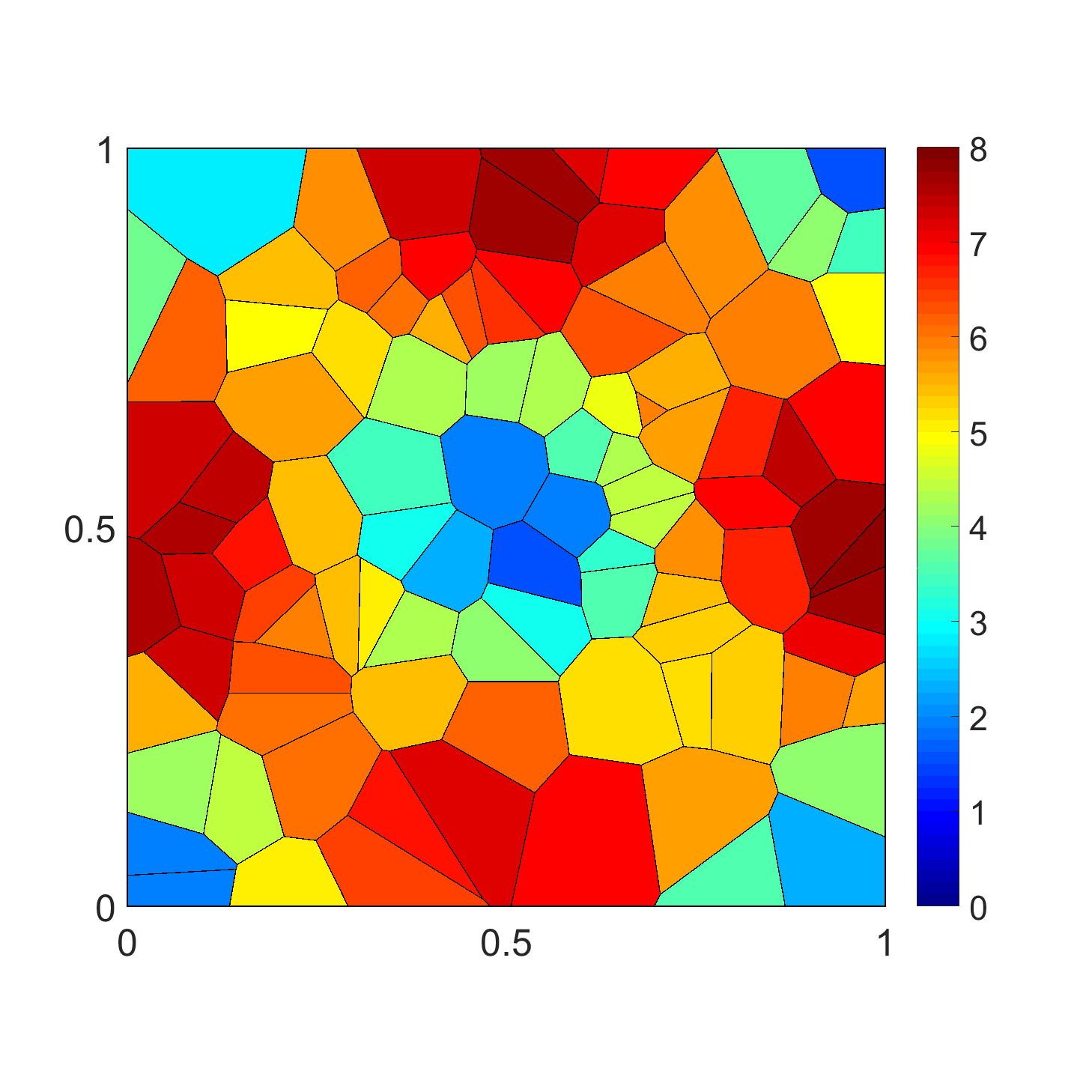}}
	\subfigure[]{\includegraphics[height=0.3\textwidth,trim = 0mm 20mm 35mm 20mm, clip]{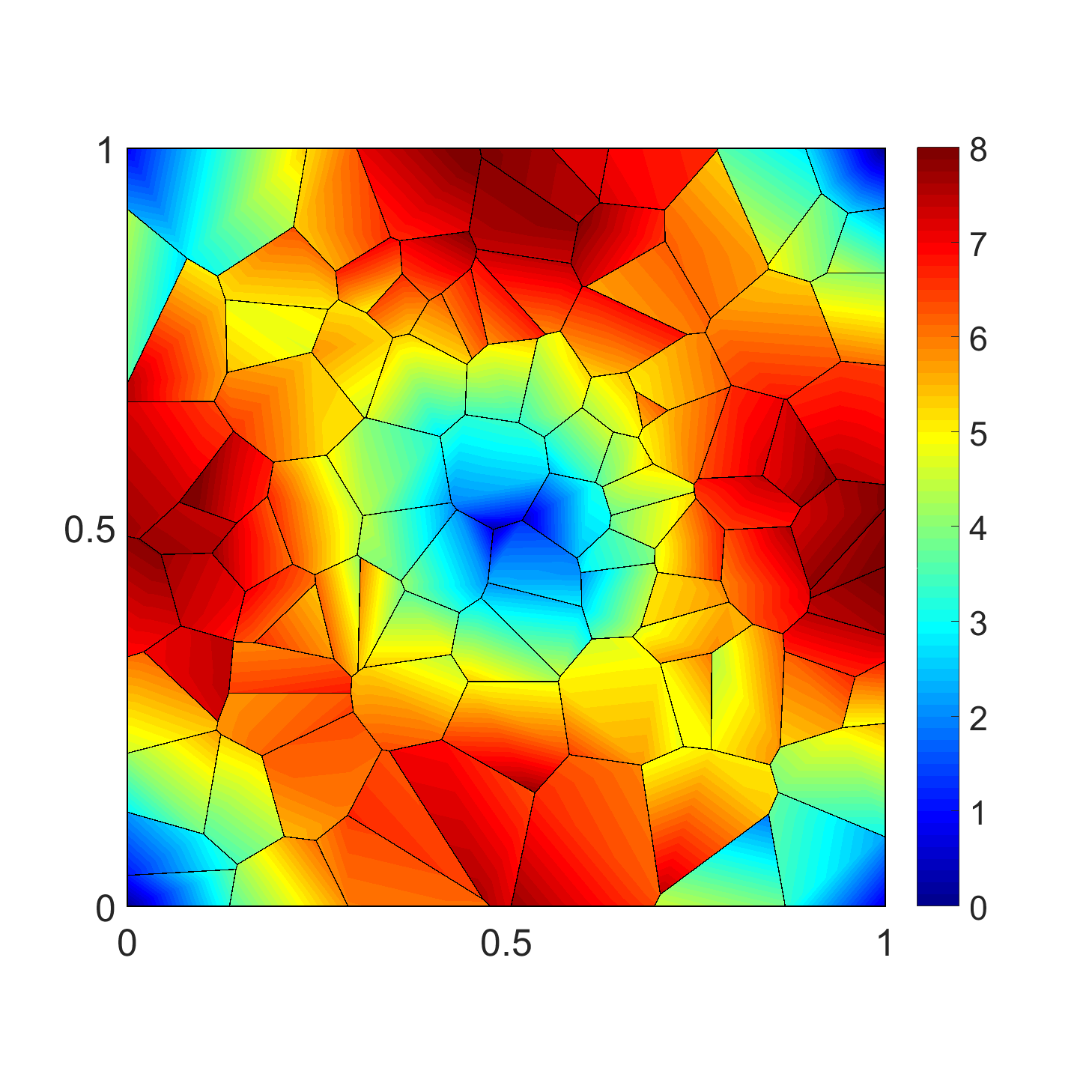}}
	\subfigure[]{\includegraphics[height=0.3\textwidth,trim = 0mm 20mm 10mm 20mm, clip]{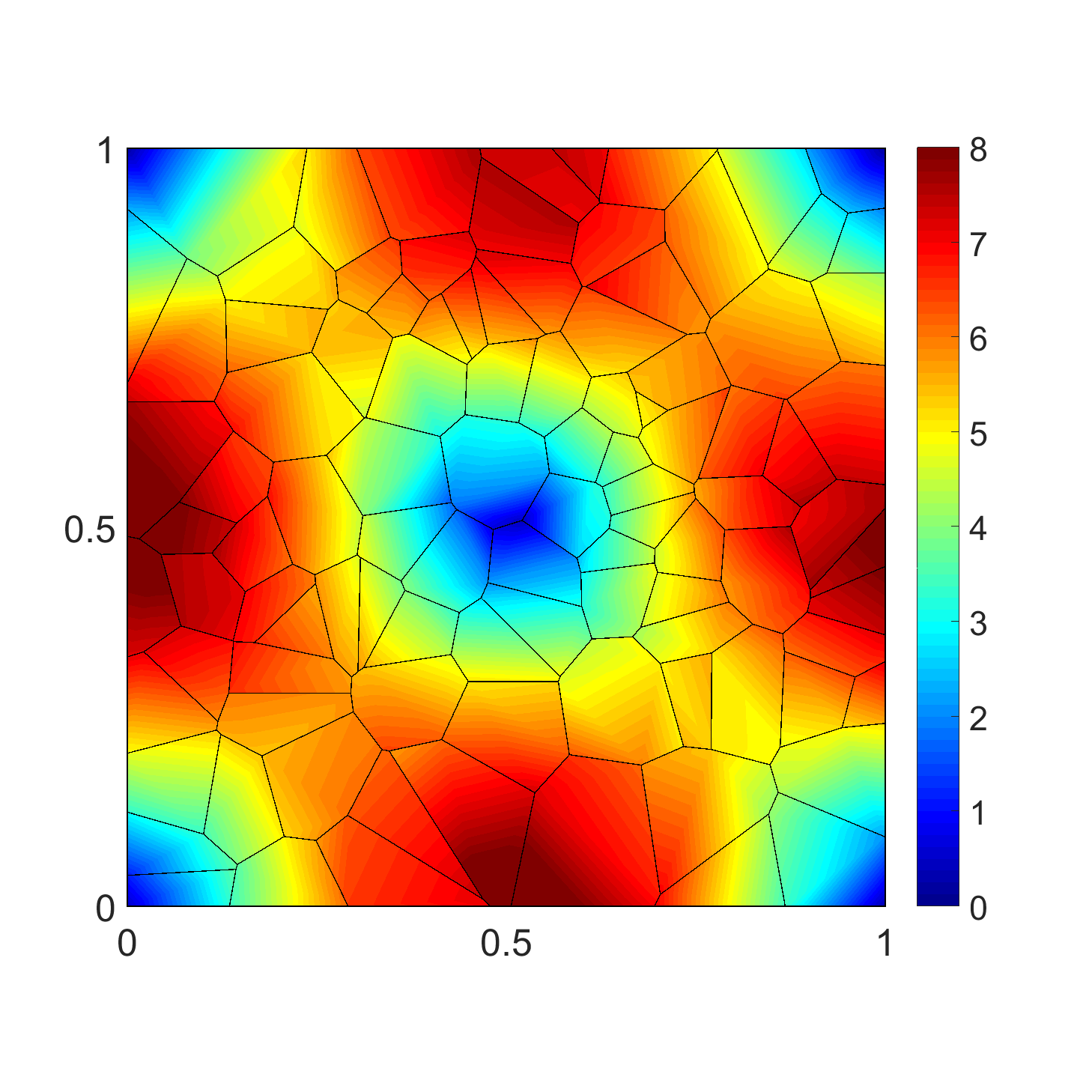}}\\
	\caption{Von Mises equivalent stress distributions for \textit{Test b}: (a) Hex (S) \textit{VEM}, (b) Hex (S) \textit{RCP0}, (c) Hex (S) \textit{RCP1}, (d) Poly (U) \textit{VEM}, (e) Poly (U) \textit{RCP0} and (f) Poly (U) \textit{RCP1}.}
	\label{fig:misesHexaTestB}
\end{figure}

Finally, a comparison between the three techniques in terms of von Mises equivalent stress distributions obtained using Hex (S) and Poly (U) meshes is reported for \textit{Test a} and \textit{Test b}. The exact solutions are shown in Fig. \ref{fig:misesAnaSol} while numerical results are reported in Figs. \ref{fig:misesHexaTestA} and \ref{fig:misesHexaTestB}, respectively. The improved compliance with the analytical solution, even for relatively coarse meshes, can be clearly observed, confirming again the effectiveness of the proposed approach.

\section{Conclusions}
\label{sec:con}

In the framework of displacement-based virtual element methods, the standard approach to compute stresses is via the constitutive law, using the available approximated strains. When general polygons are considered, such a procedure does not fully exploit the adopted degrees of freedom, thus leading to a quite poor stress field.

In the present paper, the beneficial effect of RCP used in connection with VEM discretisation schemes has been presented. The RCP formulation allows to compute accurate patchwise equilibrated stress fields starting from approximated boundary displacements and from the knowledge of the applied loads. Such a formulation appears to be well-suited for the recovery of stresses in the framework of displacement-based virtual elements, in which such quantities are known only on the elements edges.
In its simplest version, the RCP is applied elementwise, so representing an efficient alternative to the standard stress recovery. It has been also shown that further improved effects might be obtained by considering patches of elements, in agreement with the results obtained in the context of finite elements schemes.
We finally remark that RCP can be very easily integrated in existing codes without introducing significant computational costs.

\bibliography{VEM}

\end{document}